\documentclass[reqno]{amsproc}
\usepackage[latin1]{inputenc}
\usepackage{amssymb}
\usepackage{hyperref}
\usepackage{amsmath}
\usepackage{latexsym}
\usepackage{cite}
\usepackage{amssymb}
\usepackage{amsfonts}
\usepackage{amscd}
\usepackage{xcolor}
\usepackage{amstext,amsmath,amssymb,amsfonts}
\newtheorem{theorem}{Theorem}

\newtheorem{corollary}[theorem]{Corollary}

\newtheorem{definition}[theorem]{Definition}

\newtheorem{proposition}[theorem]{Proposition}
\newtheorem{remark}[theorem]{Remark}

\newcommand{\K}{\mathbb {K}}

\newcommand{\A}{\mathcal{A}}

\newcommand{\beq}{\begin{eqnarray}}
\newcommand{\eeq}{\end{eqnarray}}
\newcommand{\beqs}{\begin{eqnarray*}}
\newcommand{\eeqs}{\end{eqnarray*}}
\newcommand{\bpro}{\begin{pro}}
\newcommand{\epro}{\end{pro}}
\newcommand{\blem}{\begin{lem}}
\newcommand{\elem}{\end{lem}}
\newcommand{\bdfn}{\begin{dfn}}
\newcommand{\edfn}{\end{dfn}}
\newcommand{\bcor}{\begin{cor}}
\newcommand{\ecor}{\end{cor}}
\newcommand{\bthm}{\begin{thm}}
\newcommand{\ethm}{\end{thm}}
\newcommand{\bex}{\begin{ex}}
\newcommand{\eex}{\end{ex}}
\newcommand{\brmk}{\begin{rmk}}
\newcommand{\ermk}{\end{rmk}}
\newcommand{\bpr}{\begin{pr}}
\newcommand{\epr}{\end{pr}}
\newcommand{\benum}{\begin{enumerate}} 
\newcommand{\eenum}{\end{enumerate}}
\newcommand{\bitem}{\begin{itemize}}
\newcommand{\eitem}{\end{itemize}}

\chardef\bslash=`\\
\numberwithin{equation}{section}
\numberwithin{table}{section}
\numberwithin{theorem}{section}
\DeclareMathOperator{\id}{id}

\title[Double constructions of Heisenberg Frobenius algebras...]{Double constructions of Heisenberg Frobenius algebras and  Connes cocycles, and solutions of the three-dimensional associative Yang-Baxter equation \footnote{Preprint: ICMPA-MPA/2017/05 } }

\author{Mahouton Norbert Hounkonnou$^\ast$}
\address[$\ast$]{University of Abomey-Calavi,
International Chair in Mathematical Physics and Applications,
ICMPA-UNESCO Chair, 072 BP 50, Cotonou, Rep. of Benin}
\email{norbert.hounkonnou@cipma.uac.bj, with copy to hounkonnou@yahoo.fr}
\author{Gb\^ev\`ewou Damien  Houndedji$^\dagger$}
\address[$\dagger$]{University of Abomey-Calavi,
International Chair in Mathematical Physics and Applications,
ICMPA-UNESCO Chair, 072 BP 50, Cotonou, Rep. of Benin}
\email{ houndedjid@gmail.com}
\begin{document}
\maketitle

\today
 
\bigskip
\begin{abstract}
We consider the three-dimensional associative algebra $\mathcal{H}$ consisting of the 
$3\times 3$ strictly upper triangular matrices whose the commutator is the Heisenberg Lie 
algebra. We determine the solutions of the Yang-Baxter associative equation in 
$\mathcal{H}$. 
For  the  antisymmetric solutions, the corresponding bialgebraic structures, 
double constructions of Frobenius algebras and properties are given explicitly.
Besides, we determine some related compatible dendriform algebras and 
solutions of their $D-$equations. Using symmetric solutions of these equations, we build 
the double constructions of related Connes cocycles. Finally,  we compute solutions of the 
three-dimensional non decomposable associative Yang-Baxter equation and build the double constructions  of
associated Frobenius algebras.
\\
{
{\bf Keywords.} 
Associative algebra, bialgebra, Heisenberg algebra, Frobenius algebra, Connes cocycle, dendriform algebra, $D-$equation, Yang-Baxter equation, $ \mathcal{O} $-operator.}\\
{\bf  MSC2010.}  16T25, 05C25, 16S99, 16Z05.
\end{abstract}
\section{Introduction}
A  Frobenius algebra  is an associative
algebra equipped with a  non-degenerate invariant bilinear form. This type of algebras plays 
an important role in different areas of   mathematics and 
physics,  such as
statistical models over two-dimensional graphs \cite{[MB]} and topological quantum field theory
\cite{[JK]}. On the other hand, an antisymmetric bilinear form on an associative algebra $\mathcal A$ is an antisymmetric bilinear form on $\mathcal A$ which is a $1-$cocycle, or 
  Connes cocycle, for the Hochschid cohomology. 

 In \cite{[C.Bai5]},  C. Bai described associative analogs of Drinfeld's double constructions for Frobenius algebras and for associative algebras equipped with  non-degenerate Connes
cocycles. 
We note that there are two different types of constructions involved:
\begin{enumerate}
\item[i)] the Drinfeld's double type constructions, from a Frobenius algebra or
from an associative algebra equipped with a Connes cocyle, and
\item[ii)] the Frobenius algebra obtained from an anti-symmetric solution of the
associative Yang-Baxter equation and the non-degenerate Connes cocycle
obtained from a symmetric solution of a D-equation.
\end{enumerate} 

 The associative
Yang-Baxter equation (AYBE) and dendriform algebras appeared more early
in various  fields. Joni and Rota \cite{[Joni]} introduced the infinitesimal bialgebras to 
elaborate an algebraic framework for the calculus of divided differences. Related examples 
were provided by Aguiar in \cite{[M1]}.  
Dendriform algebras
were introduced by Loday \cite{[N8]}, inspired  by the
algebraic K-theory. Then, these structures were studied quite extensively with
connections to several areas in mathematics and physics, like
operads \cite{[N10]}, homology \cite{[N4]}, \cite{[N5]},
arithmetics \cite{[N9]}  and quantum field theory \cite{[N3]}.

In general, producing solutions of the Yang-Baxter equation (YBE) is a very hard problem, but 
very interesting and useful. Indeed, the AYBE plays an 
important role in the study of
non-abelian (Trace-) Poisson structures and non-commutative
integrable systems (see \cite{[ORS1]}, \cite{[ORS2]}, \cite{[ORS3]}, \cite{[OS]}). In 
particular, the
solutions of AYBE give the construction
of double Poisson brackets \cite{[vdB]}. 
Therefore, succeeding in providing  the
solutions of associative Yang-Baxter equations may give a good framework 
for a better  description of the above mentioned mathematical constructions and related 
physical phenomena.

In mathematics (\cite{[Binz]}, \cite{[Howe]}), the Heisenberg group, named after Werner Heisenberg, is the group of upper 
triangular matrices of the form $
\left(
\begin{array}{ccc}
1 & a & c \\
0 & 1 & b \\
0 & 0 & 1 \\  
\end{array}
\right)
$ under the operation of matrix multiplication. Elements $a, \,b$ and $c$ can be taken from any 
commutative ring with identity, often taken to be the ring of real numbers (resulting in 
the "continuous Heisenberg group") or the ring of integers (resulting in the "discrete 
Heisenberg group").
The continuous Heisenberg group arises in the description of one-dimensional quantum 
mechanical systems. More generally, one can consider Heisenberg groups associated to n-dimensional systems, and most generally, to any symplectic vector space.

In this work, we build the antisymmetric infinitesimal Heisenberg bialgebras, the double constructions of Heisenberg Frobenius algebras and Connes cocycles and determine solutions of the three-dimensional  AYBE. We start from
the vector space underlying the  Heisenberg algebra $\mathcal{H}$ generated  by the 
following matrices: \beqs
\left(
\begin{array}{ccc}
0 & 1 & 0 \\
0 & 0 & 0 \\
0 & 0 & 0 \\  
\end{array}
\right), 
\left(
\begin{array}{ccc}
0 & 0 & 1 \\
0 & 0 & 0 \\
0 & 0 & 0 \\  
\end{array}
\right),
\left(
\begin{array}{ccc}
0 & 0 & 0 \\
0 & 0 & 1 \\
0 & 0 & 0 \\  
\end{array}
\right).
\eeqs

 The paper is organized as follows. In Section 2, we give some basic notions and results on the double constructions of Frobenius algebras and Connes cocycles. In Section 3, we consider the three-dimensional associative algebra $\mathcal{H}$ consisting of the $3\times 3$ strictly upper triangular matrices whose commutator is the Heisenberg Lie algebra. We compute solutions of the Yang-Baxter associative equation in $\mathcal{H}$. For the antisymmetric solutions, the corresponding bialgebraic structures, double constructions of Frobenius algebras and properties are given explicitly. 
In section 4, we determine some related compatible dendriform algebras and solutions of their $D-$equations. Using symmetric solutions of these equations, we build the double constructions of related Connes cocycles. In section 5, given $ 3-$dimensional non decomposable complex associative algebras, we explicitly solve associative Yang-Baxter equations and use skew-symmetric solutions to perform double constructions of Frobenius algebras. Finally, in section 6, we end with some  concluding remarks.
\section{Preliminaries}
In this section, we give a quick overview of main definitions and results that we need in the sequel. For a detailed description, we refer to \cite{[C.Bai5]}. See also  \cite{[M1]},  \cite{[N14]}-\cite{[N15]}, \cite{[N8]} and  \cite{[N6]} and the references therein.
\begin{definition} \label{d1}
A bilinear form $ \mathcal{B}(\cdot, \cdot ) $ on an associative algebra $ \mathcal{A} $ is \textbf{invariant} if
\beqs
\mathcal{B}(xy,z) = \mathcal{B}(x, yz) \mbox { for all } x, y, z \in  \mathcal{A} .
\eeqs
\end{definition}

\begin{definition} \label{d2}
An antisymmetric bilinear form $ \omega(\cdot, \cdot)  $ on an associative algebra $ \mathcal{A} $ is a \textbf{cyclic 1-cocycle in the sense of Connes} if
\beq \label{eq4}
\omega(xy,z) + \omega(yz,x) + \omega(zx,y) = 0 \mbox { for all } x, y, z \in \mathcal{A}.
\eeq
For simplicity,
 $ \omega $ is called a \textbf{Connes cocycle}.
\end{definition}
In the sequel, we use the notation  $\mathcal{A}^{\ast}$ for the dual vector space $\mathcal{A}^{\ast}:= \mbox{Hom}_k(\mathcal{A}, k),$ where $k$ is the base field.

\begin{definition}  \label{d5}
We call $(\mathcal{A}, \mathcal{B})$ a   double
construction of a (symmetric) Frobenius algebra associated to
$\mathcal{A}_1$ and ${\mathcal A}_1^*$ if it satisfies the
conditions \benum \item[(1)] $ \mathcal{A} = \mathcal{A}_{1}
\oplus \mathcal{A}^{\ast}_{1} $ as the direct sum of vector
spaces; \item[(2)] $ \mathcal{A}_{1} $ and $
\mathcal{A}^{\ast}_{1} $ are associative subalgebras of $
\mathcal{A} $; \item[(3)] $ \mathcal{B} $ is the natural symmetric
bilinear form on $ \mathcal{A}_{1} \oplus \mathcal{A}^{\ast}_{1} $
given by \beq \label{eq2}
 \mathcal{B}(x + a^{\ast}, y + b^{\ast}) = \langle x, b^{\ast} \rangle +  \langle a^{\ast}, y \rangle \mbox { for all } x, y \in \mathcal{A}_{1}, a^{\ast}, b^{\ast} \in \mathcal{A}^{\ast}_{1},
\eeq
where $ \langle  , \rangle $ is the natural pair between the vector space $ \mathcal{A}_{1} $ and its dual space  $ \mathcal{A}^{\ast}_{1} $.
\eenum

We call $ (\mathcal{A}, \omega) $ 
a Connes cocycle associated to
$\mathcal{A}_1$ and ${\mathcal A}_1^*$ if it satisfies the
conditions $(1), (2)$ and \benum \item[(4)] $ \omega $ is the
natural antisymmetric bilinear form on $ \mathcal{A}_{1} \oplus
\mathcal{A}^{\ast}_{1} $ given by \beq \label{eq7} \omega(x +
a^{\ast}, y +  b^{\ast}) = -\langle x, b^{\ast} \rangle + \langle
a^{\ast}, y \rangle \mbox { for all } x, y \in \mathcal{A}_{1},
a^{\ast}, b^{\ast} \in \mathcal{A}^{\ast}_{1}, \eeq and $ \omega $
is a Connes cocycle on $ \mathcal{A} $. \eenum
\end{definition}

\begin{definition}
Let $  \mathcal{A} $ be an associative algebra and let V be a vector space. Let $ l, r ,  \mathcal{A} \rightarrow gl(V) $ be two linear maps. $ V $ (or the pair $(l,r)$, or $ (l, r, V) $) is called a \textit{bimodule} of $  \mathcal{A} $ if 
\begin{center}
$ l(xy)v = l(x)l(y)v, r(xy)v = r(y)r(x)v, l(x)r(y)v = r(y)l(x)v  $
\end{center}
for all $ x, y \in  \mathcal{A}, v \in V $.
\end{definition}
\begin{proposition} \label{pro1}
$ (l, r, V) $ is a bimodule of an associative algebra $  \mathcal{A} $ if and only if the direct sum $  \mathcal{A} \oplus V $ of vector spaces is turned into an associative algebra  by defining multiplication in $  \mathcal{A} \oplus V $ by 
\beqs
(x_{1} + v_{1}) \ast (x_{2} + v_{2}) = x_{1} \cdot x_{2} + (l(x_{1})v_{2} + r(x_{2})v_{1} )
\eeqs
for all $ x_{1}, x_{2} \in  \mathcal{A}, v_{1}, v_{2} \in V  $. We denote such an associative algebra $ ( \mathcal{A} \oplus V, \ast) $ by $\mathcal{A} \ltimes_{l, r} V $ or simply $\mathcal{A} \ltimes V $.
\end{proposition}

Let us now  give some notations useful in the sequel. Let $ \mathcal{A} $ be an associative algebra.

 Consider the representations of the  left $L$ and right $R$ multiplication operations defined as:
 \begin{eqnarray}
 L: \A & \longrightarrow & \mathfrak{gl}(\A)  \cr
  x  & \longmapsto & L_x:
  \begin{array}{ccc}
 \A &\longrightarrow & \A \cr 
  y & \longmapsto & x \cdot y, 
   \end{array}
\end{eqnarray}

\begin{eqnarray}
    R: \A & \longrightarrow & \mathfrak{gl}(\A)  \cr
     x  & \longmapsto & R_x:
     \begin{array}{ccc}
    \A &\longrightarrow & \A \cr 
     y & \longmapsto & y \cdot x,
      \end{array}
 \end{eqnarray}
 The dual maps $L^{*}, R^{*}$  of the linear maps $L, R,$ are defined, respectively, as: 
$\displaystyle L^{*}, R^{*}: \A \rightarrow \mathfrak{gl}(\A^{*})$
 such that:
\beq\label{dual1}
 L^*: \A & \longrightarrow & \mathfrak{gl}(\A^*)  \cr
  x  & \longmapsto & L^*_x:
      \begin{array}{llll}
 \A^* &\longrightarrow & \A^* \\ 
  u^* & \longmapsto & L^*_x u^*: 
      \begin{array}{llll}
\A  &\longrightarrow&  \K \cr
v  &\longmapsto&  L^{*}_x(u^{*})(v) 
 :=   u^{*}(x\cdot v), 
      \end{array}
     \end{array}
 \eeq
\beq\label{dual2}
 R^*: \A & \longrightarrow & \mathfrak{gl}(\A^*)  \cr
  x  & \longmapsto & R^*_x:
      \begin{array}{llll}
 \A^* &\longrightarrow & \A^* \\ 
  u^* & \longmapsto & R^*_x u^*: 
      \begin{array}{lllll}
\A  &\longrightarrow&  \K \cr
v  &\longmapsto&  R^{*}_x(u^{*})v 
 :=  u^{*}(v\cdot x), 
      \end{array}
     \end{array}
\eeq
for all $x, v \in \A, u^{*} \in \A^{*},$ where $\A^{*}$ is the dual space of $\A.$

 Let $ \sigma : \mathcal{A}\otimes \mathcal{A} \rightarrow \mathcal{A}\otimes \mathcal{A}  $ be the exchange operator defined as
\beqs
\sigma(x \otimes y) = y \otimes x,
\eeqs for all $ x, y \in \mathcal{A} $.
  A solution of the associative Yang-Baxter equation  (AYBE) in the associative algebra $\mathcal{A}$ is an element  $ r = \sum_{i} x_{i} \otimes y_{i} \in \mathcal{A} \otimes \mathcal{A}$
satisfying the equation \cite{[C.Bai5]}
\beq \label{AYEB}
r_{12}r_{13} + r_{13}r_{23} - r_{23}r_{12} = 0,
\eeq
where 
\beqs
r_{12}r_{13} = \sum_{i,j} x_{i}x_{j} \otimes y_{i} \otimes y_{j},
\eeqs
\beqs
r_{13}r_{23} = \sum_{i,j} x_{i} \otimes x_{j} \otimes y_{i}y_{j},
\eeqs
\beqs
r_{23}r_{12} = \sum_{i,j} x_{j} \otimes x_{i}y_{j} \otimes y_{i}.
\eeqs

\begin{definition}
Let $V_{1} $, $ V_{2} $ be two vector spaces. For a linear map $ \phi : V_{1} \rightarrow V_{2} $, we denote the dual (linear) map by $ \phi^{\ast} : V^{\ast}_{2} \rightarrow V^{\ast}_{1} $ given by 
\beqs
\langle v, \phi^{\ast}(u^{\ast})\rangle = \langle \phi(v), u^{\ast} \rangle
\eeqs 
for all  $ v \in V_{1} $, $ u^{\ast} \in V^{\ast}_{2} $.
\end{definition}

\begin{theorem}\cite{[C.Bai5]} \label{pro2.2.2}
Let $ (\mathcal{A}, \cdot) $ be an associative algebra. Suppose that there is an associative algebra 
 structure $ "\circ" $ on its dual space $ \mathcal{A}^{\ast} $. Then,  $ (\mathcal{A}, \mathcal{A}^{\ast},
  R^{\ast}, L^{\ast}, R^{\ast}_{\circ},  L^{\ast}_{\circ}) $ is a matched pair of associative algebras 
  if and only if for any $ x \in \mathcal{A}$ and $ a^{\ast}, b^{\ast} \in \mathcal{A}^{\ast} $,
\beq \label{eq16}
R^{\ast}(x)(a^{\ast} \circ b^{\ast}) = R^{\ast}(L^{\ast}_{\circ}(a^{\ast})x)b^{\ast} + 
(R^{\ast}(x)a^{\ast})\circ b^{\ast}, 
\eeq
\beq \label{eq17}
 R^{\ast}(R^{\ast}_{\circ}(a^{\ast})x)b^{\ast} + L^{\ast}(x)a^{\ast}\circ b^{\ast}=  L^{\ast}(L^{\ast}_{\circ}(b^{\ast})x)a^{\ast} + a^{\ast}\circ (R^{\ast}(x)b^{\ast}).
\eeq
\end{theorem}

\begin{definition} An infinitesimal bialgebra is a triple $(\mathcal{A}, \cdot, \Delta),$ where $(\mathcal{A}, \cdot)$ is an associative algebra, and $(\mathcal{A}, \Delta)$ is a coassociative coalgebra and

\begin{eqnarray}
\Delta(ab)= \sum ab_1\otimes b_2 + \sum a_1\otimes a_2 b \,\,\,\forall \,a,\, b \in \mathcal{A},
\end{eqnarray}
or, equivalently, in terms of left and right multiplication operators,
\begin{eqnarray}
\Delta(ab)= \left(L(a)\otimes I) \Delta(b) + (I\otimes R(b)\right)\Delta(a),
\end{eqnarray}
where $I$ is the identity operator, and we denote $\Delta(a)= \sum a_1\otimes a_2$, for any element $a \in \mathcal{A}.$
\end{definition}

\begin{definition} \label{aiba} 
 Let $ \mathcal{A} $ be a finite dimensional  associative algebra over the field $k$. An \textbf{antisymmetric infinitesimal
 bialgebra} structure on $ \mathcal{A} $ is a linear map $ \Delta: \mathcal{A} \rightarrow
\mathcal{A} \otimes \mathcal{A} $ such that
\begin{enumerate}
\item $  \Delta^{\ast} : \mathcal{A}^{\ast} \otimes \mathcal{A}^{\ast}  \rightarrow \mathcal{A}^{\ast}  $ defines an associative algebra structure on $ \mathcal{A}^{\ast} $;
\item $ \Delta $ satisfies the following equations:
\begin{eqnarray}
&&\Delta (x\cdot y) = (I\otimes L(x))\Delta(y) + (R(y) \otimes I)\Delta(x),\\
&&
(L(y) \otimes I - I\otimes R(y))\Delta(x)
+ \sigma [(L(x)\otimes I - I\otimes R(x))\Delta(y)] = 0,
\end{eqnarray}
for all $ x, y \in \mathcal{A} $.
\end{enumerate}
\end{definition}
We denote this bialgebra structure by $ (\mathcal{A}, \Delta) $ or simply by $(\mathcal{A},
 \mathcal{A}^{\ast}).$ This last notation, i. e. $(\mathcal{A},
 \mathcal{A}^{\ast}),$ only used in this work as a matter of  simplification, has not to be mistaken with the associative algebra $\mathcal{A}\oplus  \mathcal{A}^{\ast},$ which is the direct sum of  $\mathcal{A}$ and $ \mathcal{A}^{\ast}.$
 
\begin{theorem}\cite{[C.Bai5]} \label{c2}
Let $(\mathcal{A}, \cdot)$ and $(\mathcal{A}^{\ast}, \circ)$ be two associative algebras. Then, the following conditions are equivalent:
\begin{enumerate}
\item there is a double construction of a Frobenius algebra associated with  $(\mathcal{A}, \cdot) $ and $(\mathcal{A}^{\ast}, \circ)$;
\item $(\mathcal{A}, \mathcal{A}^{\ast}, R^{\ast}_{\cdot}, L^{\ast}_{\cdot}, R^{\ast}_{\circ},  L^{\ast}_{\circ}) $ is a matched pair of associative algebras; 
\item $(\mathcal{A}, \mathcal{A}^{\ast})$ is an antisymmetric infinitesimal bialgebra.
\end{enumerate}
\end{theorem}

\begin{corollary}\label{c2}
Let $ \mathcal{A} $ be an associative algebra, and $ r \in \mathcal{A} \otimes \mathcal{A} $ be an antisymmetric solution of the AYBE (\ref{AYEB}) in $\mathcal{A}. $
The map $ \Delta $ defined by

\begin{eqnarray}
\Delta(x) = (I\otimes L(x) - R(x)\otimes I)r \mbox { for all }  x \in \mathcal{A}
\end{eqnarray}

 induces an associative algebra structure on $ \mathcal{A}^{\ast}$ and
 $ (\mathcal{A}, \mathcal{A}^{\ast}) $ is an antisymmetric infinitesimal bialgebra in the sense of {\rm Definition} \ref{aiba}. 
\end{corollary}

\begin{proposition}\cite{[C.Bai5]}  \label{pro2.4.4}
Let $(\mathcal{A},\cdot)$ be an associative algebra and let $ r
\in \mathcal{A} \otimes \mathcal{A} $ be an antisymmetric solution
of the AYBE in $ \mathcal{A} $. Then, the corresponding double construction of Frobenius algebra
$(\mathcal{AD}(\mathcal{A}),\ast) $ associated to $\mathcal{A}$
and $\mathcal A^*$
 is given from the
product in $ \mathcal{A} $ as follows: \beq \label{eq27} a^{\ast}
\ast b^{\ast} = a^{\ast} \circ b^{\ast} =
R^{\ast}(r(a^{\ast}))b^{\ast} + L^{\ast}(r(b^{\ast})) a^{\ast},
\eeq \beq \label{eq28} x \ast a^{\ast} = x \cdot r(a^{\ast}) -
r(R^{\ast}(x)a^{\ast}) + R^{\ast}(x)a^{\ast}, \eeq \beq
\label{eq29} a^{\ast} \ast x = r(a^{\ast}) \cdot x -
r(L^{\ast}(x)a^{\ast}) + L^{\ast}(x)a^{\ast}, \eeq for any $ x \in
\mathcal{A}, a^{\ast}, b^{\ast} \in \mathcal{A}^{\ast} $.
\end{proposition}

\begin{definition}
Let $ \mathcal{A} $ be a vector space
with two bilinear products denoted by $
\prec $ and $ \succ $. Then $ (\mathcal{A}, \prec, \succ) $ is
called a \textbf{dendriform algebra} if, for any $ x, y, z \in
\mathcal{A} $,
\begin{eqnarray*}
(x \prec y) \prec z &=& x \prec (y \ast z), \cr
(x \succ y) \prec z &=& x \succ (y \prec z), \cr
x \succ (y \succ z) &=& ( x \ast y) \succ z,
\end{eqnarray*}
where $ x \ast y = x \prec y + x \succ y $.
\end{definition}

Let $ (\mathcal{A}, \prec, \succ) $ be a dendriform algebra. For any $ x \in \mathcal{A} $, let $ L_{\prec}(x),  R_{\prec}(x) $ and $ L_{\succ}(x),$ \\$ R_{\succ}(x) $ denote the left and right multiplication operators of $(\mathcal{A}, \prec)$ and $(\mathcal{A}, \succ)$, respectively:
\begin{eqnarray}\nonumber
L_{\prec}(x) y = x \prec y, 
R_{\prec}(x) y = y \prec x, L_{\succ}(x) y = x \succ y, R_{\succ}(x) y = y \succ x, 
\end{eqnarray}
for all $ x, y \in \mathcal{A} $. Moreover, let $ L_{\prec}, R_{\prec}, L_{\succ}, R_{\succ} : \mathcal{A} \rightarrow gl(\mathcal{A})$ be four linear maps with $ x \mapsto L_{\prec}(x),  x \mapsto R_{\prec}(x),  x \mapsto L_{\succ}(x), $ and $  x \mapsto R_{\succ}(x)$, respectively. It is known that the product given by \cite{[N8]} 
\begin{eqnarray} \label{eq32}
x \ast y = x \prec y + x \succ y, \mbox { for all } x, y \in \mathcal{A},
\end{eqnarray}
is associative. We call $ (\mathcal{A}, \ast) $ the associated associative algebra of $ (\mathcal{A}, \prec, \succ).$ Conversely $ (\mathcal{A}, \prec, \succ)$ is called a compatible dendriform algebra structure on the associative algebra $ (\mathcal{A}, \ast) $.

\begin{theorem}\cite{[C.Bai5]} \label{theo4.1.1}
Let $(\mathcal{A}, \ast)$ be an associative algebra and let $ \omega $ be a non-degenerate Connes cocycle. 
Then,  there exists a compatible dendriform algebra structure $ \succ, \prec $ on $  \mathcal{A} $ given by 
\beq \label{eq53}
\omega(x \succ y, z) = \omega(y, z \ast x), \ \ \ \omega(x \succ y, z) = \omega(x, y \ast z) \mbox { for all } x, y \in \mathcal{A}.
\eeq 
\end{theorem}

\begin{corollary} Let $(T(\mathcal{A}) = \mathcal{A} \bowtie \mathcal{A}^{\ast}, \omega)$ be a double construction of the Connes cocycle. Then, 
 there exists a compatible dendriform algebra structure $ \succ, \prec $ on $ T(\mathcal{A}) $ 
defined by the equation (\ref{eq53}).
  Moreover, $ \mathcal{A} $ and $ \mathcal{A}^{\ast} $, endowed with this product,
 are dendriform subalgebras. 
\end{corollary}

\begin{definition} 
Let $ \mathcal{A} $ be a vector space. A \textbf{dendriform D-bialgebra} structure on $ \mathcal{A} $ is a set of linear maps
$ (\Delta_{\prec}, \Delta_{\succ}, \beta_{\prec}, \beta_{\succ}) $  given by
 \ \ $ \Delta_{\prec}, \Delta_{\succ} : \mathcal{A} \rightarrow \mathcal{A} \otimes \mathcal{A} $, \ \
$ \beta_{\prec},  \beta_{\succ} : \mathcal{A}^{\ast} \rightarrow \mathcal{A}^{\ast} \otimes \mathcal{A}^{\ast} $, such that
\benum
\item[(a)] $ (\Delta^{\ast}_{\prec}, \Delta^{\ast}_{\succ}) : \mathcal{A}^{\ast} \otimes \mathcal{A}^{\ast} \rightarrow \mathcal{A}^{\ast} $
 defines a dendriform algebra structure $(\succ_{\mathcal{A}^{\ast}}, \prec_{\mathcal{A}^{\ast}}) $ on $ \mathcal{A}^{\ast} $;
\item[(b)] $ (\beta^{\ast}_{\prec}, \beta^{\ast}_{\succ}) : \mathcal{A} \otimes \mathcal{A} \rightarrow A $ defines a dendriform algebra
 structure $(\succ_{\mathcal{A}}, \prec_{\mathcal{A}})$ on $  \mathcal{A} $;
\item[(c)] the following equations are satisfied
\beq \label{eq55}
\Delta_{\prec}(x \ast_{\mathcal{A}} y) = (\id \otimes L_{\prec_{\mathcal{A}}}(x))\Delta_{\prec}(y) + (R_{\mathcal{A}}(y)\otimes \id)\Delta_{\prec}(y),
\eeq
\beq \label{eq56}
\Delta_{\succ}(x \ast_{\mathcal{A}} y) = (\id \otimes L_{\prec_{\mathcal{A}}}(x))\Delta_{\succ}(y) + (R_{\prec_{A}}(y)\otimes \id)\Delta_{\succ}(y),
\eeq
\beq \label{eq57}
\beta_{\prec}(a^{\ast} \ast_{\mathcal{A}^{\ast}} b^{\ast}) = (\id \otimes L_{\prec_{\mathcal{A}^{\ast}}}(a^{\ast}))\beta_{\prec}(b^{\ast}) + (R_{\mathcal{A}^{\ast}}(b^{\ast})\otimes \id)\beta_{\prec}(a^{\ast})
\eeq
\beq \label{eq58}
\beta_{\succ}(a^{\ast} \ast_{\mathcal{A}^{\ast}} b^{\ast}) = (\id \otimes L_{\mathcal{A}^{\ast}}(a^{\ast}))\beta_{\succ}(b^{\ast}) + (R_{\prec_{\mathcal{A}^{\ast}}}(b^{\ast})\otimes \id)\beta_{\succ}(a^{\ast}) ,
\eeq
\beq \label{eq59}
(L_{\mathcal{A}}(x)\otimes \id - \id \otimes R_{\prec_{\mathcal{A}}}(x))\Delta_{\prec}(y) + \sigma[(L_{\succ_{\mathcal{A}}}(y)\otimes (-\id)\otimes R_{\mathcal{A}}(y))\Delta_{\prec}(y)] = 0 ,
\eeq
\beq \label{eq60}
(L_{\mathcal{A}^{\ast}}(a^{\ast})\otimes \id - \id \otimes R_{\prec_{\mathcal{A}^{\ast}}}(a^{\ast}))\beta_{\prec}(b^{\ast}) 
+ \sigma[(L_{\succ_{\mathcal{A}^{\ast}}}(b^{\ast})\otimes(-id)\otimes R_{\mathcal{A}^{\ast}}(b^{\ast}))\beta_{\succ}(a^{\ast})] = 0,
\eeq
hold for any $ x, y \in \mathcal{A} $ and
$ a^{\ast}, b^{\ast} \in \mathcal{A}^{\ast} $,
where $ L_{\mathcal{A}} = L_{\succ_{\mathcal{A}}} + L_{\prec_{\mathcal{A}}}, R_{\mathcal{A}}
= R_{\succ_{\mathcal{A}}} + R_{\prec_{\mathcal{A}}}, L_{\mathcal{A}^{\ast}} = L_{\succ_{\mathcal{A}^{\ast}}}
+ L_{\prec_{\mathcal{A}^{\ast}}},  R_{\mathcal{A}^{\ast}} = R_{\succ_{\mathcal{A}^{\ast}}} + R_{\prec_{\mathcal{A}^{\ast}}} $.

\eenum
We also  denote it by $ (\mathcal{A}, \mathcal{A}^{\ast}, \Delta_{\succ}, \Delta_{\prec},  \beta_{\succ}, \beta_{\prec}) $
 or simply $(\mathcal{A}, \mathcal{A}^{\ast})$.
\end{definition}

\begin{theorem}\cite{[C.Bai5]}
Let $(A, \prec_{A}, \succ_{A})$ and $(A^{\ast}, \prec_{A^{\ast}}, \succ_{A^{\ast}} )$ be two dendriform algebras. Let $(A, \ast_{A})$ and   $(A^{\ast}, \ast_{A^{\ast}})$ be the associated associative algebras, respectively. Then, the following
 conditions are equivalent:
 \begin{enumerate}
 \item there is a double construction of the Connes cocycle associated with
 $(A, \ast_{A})$ and $(A^{\ast}, \ast_{A^{\ast}})$;
\item $ (A, A^{\ast}) $ is a dendriform $ D- $bialgebra.
\end{enumerate}
\end{theorem}

\begin{corollary} \label{cor1}
Let $ (\mathcal{A}, \succ, \prec) $ be a dendriform algebra and
 $ r \in \mathcal{A} \otimes \mathcal{A} $. 
Suppose that
 $ r  $ is symmetric and $ r $ satisfies the equation
\begin{eqnarray}
\label{eqD}
r_{12} \ast r_{13} = r_{13} \prec r_{23} + r_{23} \succ r_{12}.
\end{eqnarray}
Then, the maps $  \Delta_{\succ}$ and $\Delta_{\prec} $ are  defined, respectively, by
\begin{eqnarray}
\Delta_{\succ}(x) &=& (id \otimes L(x) - R_{\prec}(x) \otimes id)r_{\succ}, \cr
\Delta_{\prec}(x) &=& (id \otimes L_{\succ}(x) - R(x) \otimes id)r_{\prec}, \forall x
\in \mathcal{A},
\end{eqnarray}
 where  $ r_{\succ} = - r$ and $ r_{\prec} = r$ induce a dendriform algebra structure on $ \mathcal{A}^{\ast} $ such that $(\mathcal{A}, \mathcal{A}^{\ast})$ is a dendriform $D$-bialgebra. Equation (\ref{eqD}) is called a \textbf{$D$-equation} in $ \mathcal{A} $.
\end{corollary}
\begin{proposition}\cite{[C.Bai5]}  \label{pro4.4.6}
Let $ (\mathcal{A}, \succ, \prec) $ be a dendriform algebra  and let $ r \in \mathcal{A} \otimes \mathcal{A} $
 be a symmetric solution of the $D$-equation in $ \mathcal{A} $. Then,  the corresponding double
construction of Connes cocycle associated to $\mathcal{A}$ and
$\mathcal A^*$ is given from the products in $ \mathcal{A} $ as
follows:
 \beq \label{eq80} 
 a^{\ast} \prec b^{\ast} &=&
-R^{\ast}_{\succ}(r(a^{\ast}))b^{\ast} +
L^{\ast}(r(b^{\ast}))a^{\ast}, \cr 
a^{\ast} \succ b^{\ast} &=&
R^{\ast}(r(a^{\ast}))b^{\ast}
-L^{\ast}_{\prec}(r(b^{\ast}))a^{\ast}, \cr
a^{\ast} \ast
b^{\ast} &=& a^{\ast} \succ b^{\ast} + a^{\ast} \prec b^{\ast} =
R^{\ast}_{\prec}(r(a^{\ast}))b^{\ast} +
L^{\ast}_{\succ}(r(b^{\ast}))a^{\ast}, \cr
x \succ a^{\ast} &=& x
\succ r(a^{\ast}) - r(R^{\ast}(x)a^{\ast}) + R^{\ast}(x)a^{\ast},
\cr 
x \prec a^{\ast} &=&  x \prec r(a^{\ast}) +
r(R^{\ast}_{\succ}(x)a^{\ast}) - R^{\ast}_{\succ}(x)a^{\ast}, \cr
x \ast a^{\ast} &=& x \ast r(a^{\ast}) -
r(R^{\ast}_{\prec}(x)a^{\ast}) + R^{\ast}_{\prec}(x)a^{\ast}, \cr
a^{\ast} \succ x &=& r(a^{\ast}) \succ x +
r(L^{\ast}_{\prec}(x)a^{\ast}) -  L^{\ast}_{\prec}(x)a^{\ast}, \cr
a^{\ast} \prec x &=& r(a^{\ast}) \prec x -  r(L^{\ast}(x)a^{\ast})
+ L^{\ast}(x)a^{\ast}, \cr
a^{\ast} \ast x &=& r(a^{\ast}) \ast x
-  r(L^{\ast}_{\succ}(x)a^{\ast}) + L^{\ast}_{\succ}(x)a^{\ast}
\eeq for any $ x \in \mathcal{A}, a^{\ast}, b^{\ast} \in
\mathcal{A}^{\ast} $.
\end{proposition}

In the sequel, unless otherwise stated, all the parameters belong to the complex field $\mathbb{C}$.

 \section{Double constructions of  Heisenberg Frobenius  algebras}
 \subsection{Heisenberg associative algebra}
\begin{definition}
Let $X, Y, Z$ be the matrices defined by:
\beqs
X =  \left(
\begin{array}{ccc}
0 & 1 & 0 \\
0 & 0 & 0 \\
0 & 0 & 0 \\  
\end{array}
\right), 
Y = \left(
\begin{array}{ccc}
0 & 0 & 1 \\
0 & 0 & 0 \\
0 & 0 & 0 \\  
\end{array}
\right),
Z = \left(
\begin{array}{ccc}
0 & 0 & 0 \\
0 & 0 & 1 \\
0 & 0 & 0 \\  
\end{array}
\right).
\eeqs
The real Lie algebra $ g $ generated by $ X, Y, Z $ is called the \textbf{Heisenberg algebra}, and the group $ G $ of real matrices 
of the form \beqs
\left(
\begin{array}{ccc}
1 & a & c \\
0 & 1 & b \\
0 & 0 & 1 \\  
\end{array}
\right)
\eeqs is called the \textbf{Heisenberg group}.
\end{definition}
Set $ \mathcal{H} = $ $ g $, $ e_{1} = X $, $ e_{2} = Y $, $ e_{3} = Z  $ and " $ \cdot $ " the  matrix product. 
\begin{proposition}
$ (\mathcal{H}, \cdot) $ is a non commutative nilpotent associative algebra. Besides,
\beq \label{eq106}
e_{i} \cdot e_{j} = \delta_{ik} \delta_{jl} e_{2}, \ \ i, j=1, 2, 3; \ \ k=1, l= 3.
\eeq
\end{proposition}
We call it the \textbf{Heisenberg associative algebra}.

\textbf{Proof: }
Let $ A = \left(
\begin{array}{ccc}
0 & a & b \\
0 & 0 & c \\
0 & 0 & 0 \\  
\end{array}
\right), B = \left(
\begin{array}{ccc}
0 & a' & b' \\
0 & 0 & c' \\
0 & 0 & 0 \\  
\end{array}
\right) \in \mathcal{H} $. We have
\beqs
A\cdot B =  \left(
\begin{array}{ccc}
0 & 0 & ac' \\
0 & 0 & 0 \\
0 & 0 & 0 \\  
\end{array}
\right) = ac' \left(
\begin{array}{ccc}
0 & 0 & 1 \\
0 & 0 & 0 \\
0 & 0 & 0 \\  
\end{array} 
\right) = ac' e_{2} \in \mathcal{H}.
\eeqs
Then $ \cdot $ is an operation in $ \mathcal{H} $. Hence,  the equation (\ref{eq106}) is obtained
 by straightforward computation of the products of matrices, $ e_{i}, e_{j}, (i, j= 1, 2, 3) $. Now, let $ C \in \mathcal{H} $.  We have
\beqs
(A\cdot B)\cdot C = 0 = A\cdot (B\cdot C). 
\eeqs
Then the product $'' \cdot''$ is associative.

$ \hfill \square $

\begin{proposition}
Let  $ (\mathcal{H}_{\alpha}, \cdot) $ be the complex algebra defined by
\beq
&&e_{i} \cdot e_{j} = \delta_{ik} \delta_{jl} e_{2}, \ \ i, j=1, 2, 3; \ \ k=1, l= 3;\cr
&&e_{m} \cdot e_{n} =\alpha \delta_{mu} \delta_{nv} e_{2}, \ \ m, n=1, 2, 3; \ \ v=1, u= 3; \alpha \in \mathbb{C}-\lbrace 1 \rbrace.
\eeq
$ (\mathcal{H}_{\alpha}, \cdot) $ is a non commutative nilpotent associative algebra.
\end{proposition}
We call it the \textbf{extended Heisenberg associative algebra}.
\begin{remark}
For $ \alpha= 1 $, $ (\mathcal{H}_{1}, \cdot) $ is a commutative nilpotent associative algebra. We call it the \textbf{commutative Heisenberg associative algebra}.
\end{remark}
\subsection{Solutions of the associative Yang-Baxter equation in $\mathcal{H}_{\alpha}$}

Let $ r $ be an element of $ \mathcal{H}_{\alpha}\otimes \mathcal{H}_{\alpha} $. Set 
\beqs
r &=& a_{11} e_{1}\otimes e_{1} + a_{12} e_{1}\otimes e_{2} + a_{13} e_{1}\otimes e_{3} + a_{21} e_{2}\otimes e_{1}\cr
&& + a_{22} e_{2}\otimes e_{2} + a_{23} e_{2}\otimes e_{3}
 + a_{31} e_{3}\otimes e_{1} + a_{32} e_{3}\otimes e_{2} + a_{33} e_{3}\otimes e_{3}. 
\eeqs

Then, the AYBE provides the  equations in Table 1:
\begin{center}
{\bf Table 1}: Associative Yang-Baxter equations in $\mathcal{H}_{\alpha}$
\begin{tabular}{c|c}
  \hline
 $E_{1}: a_{11}a_{31} + \alpha a_{11}a_{31}= 0$ & $E_{2}: a_{11}a_{32} + \alpha a_{31}a_{12} + a_{21}a_{13}= 0$ 
\\ \hline
$E_{3}: a_{11}a_{33} + \alpha a_{31}a_{13}= 0$ & $E_{4}: a_{12}a_{31} + \alpha a_{32}a_{11}= a_{11}a_{23} + \alpha a_{31}a_{21}$
\\ \hline
$E_{5}: \alpha a_{23}a_{11}= 0$ & $ E_{6}: a_{12}a_{33} + \alpha a_{32}a_{13} + a_{21}a_{23}= a_{13}a_{23} + \alpha a_{33}a_{21}$
\\ \hline
$E_{7}: a_{13}a_{31} + \alpha a_{33}a_{11}= 0$ & $E_{8}: a_{13}a_{32} + \alpha a_{33}a_{12} + a_{21}a_{33} + \alpha a_{23}a_{31} = 0$
\\ \hline
$E_{9}: a_{13}a_{33} + \alpha a_{33}a_{13}= 0$ & $E_{10}: a_{11}a_{13} + \alpha a_{13}a_{11}= 0$
\\ \hline
$E_{11}: a_{11}a_{33} + \alpha a_{13}a_{31}= 0$ & $E_{12}: a_{31}a_{13} + \alpha a_{33}a_{11}= 0$
\\ \hline
$E_{13}: a_{11}a_{23} + \alpha a_{13}a_{21}= $ & $E_{14}: a_{31}a_{23} + \alpha a_{33}a_{21}= $\cr
 $a_{12}a_{13} + \alpha a_{32}a_{11}$ & $a_{12}a_{33} + \alpha a_{32}a_{31}$
\\ \hline
$E_{15}: a_{31}a_{33} + \alpha a_{33}a_{31}= 0$ & $E_{16}: a_{12}a_{32} + \alpha a_{32}a_{12} + \alpha a_{23}a_{21}=$\cr
& $a_{12}a_{23} + \alpha a_{32}a_{21}$
\\ \hline
$ E_{17}: a_{11}a_{13} + \alpha a_{31}a_{11}= 0$ & $E_{18}: a_{11}a_{33} + \alpha a_{31}a_{31}= 0$
\\ \hline
$ E_{19}: a_{13}a_{13} + \alpha a_{33}a_{11}= 0$ & $E_{20}: a_{13}a_{33} + \alpha a_{33}a_{31}= 0$
\\ \hline
\end{tabular}
\end{center}
\begin{proposition}
 The solutions of the corresponding AYBE in $ (\mathcal{H}_{\alpha}, \cdot) $ are given in  Table 2.
\begin{proposition}
The  antisymmetric solutions of the AYBE in $ \mathcal{H}_{\alpha} $ are
\beq
r_{1}= a_{12}(e_{1}\otimes e_{2} - e_{2}\otimes e_{1}) \mbox{ and } 
r_{2}= a_{23}(e_{2}\otimes e_{3} - e_{3}\otimes e_{2}) .
\eeq
\end{proposition}

\begin{center}
{\bf Table 2}:
 Solutions of the  associative Yang-Baxter equation in $ (\mathcal{H}_{\alpha}, \cdot)$
\begin{tabular}{c c}
\hline
Associative & Solutions of the AYBE \cr
 algebra $ \mathcal{A} $   &
 \\ \hline
$ \mathcal{H}_{1} $: 
 $ e_{1}\cdot e_{3}= e_{2} $ & $ \left(
\begin{array}{ccc}
0 & a_{12} & 0 \\
a_{21} & a_{22} & 0 \\
0 & 0 & 0 \\  
\end{array}
\right)
$; $ \left(
\begin{array}{ccc}
0 & 0 & 0 \\
0 & a_{22} & a_{23} \\
0 & a_{32} & 0 \\  
\end{array}
\right) a_{23}\neq 0
$               \cr       
$ e_{3}\cdot e_{1}= e_{2} $ & $ \left(
\begin{array}{ccc}
0 & 0 & 0 \\
0 & a_{22} & a_{23} \\
0 & a_{32} & a_{33} \\  
\end{array}
\right) a_{33}\neq 0;
$ 
  $ \left(
\begin{array}{ccc}
a_{11} & a_{12} & 0 \\
a_{21} & a_{22} & 0 \\
0 & 0 & 0 \\  
\end{array}
\right) a_{11}\neq 0
 $ \cr 
 & $ \left(
\begin{array}{ccc}
0 & a_{12} & 0 \\
a_{21} & a_{22} & 0 \\
0 & a_{32} & 0 \\  
\end{array}
\right)  a_{21}= 2a_{12}, a_{32}\neq 0
 $ \cr
 & $ \left(
\begin{array}{ccc}
0 & a_{12} & 0 \\
0 & a_{22} & a_{23} \\
0 & a_{32} & 0 \\  
\end{array}
\right) a_{23}= 2a_{32} ,a_{23}, a_{12}\neq 0$
\\ \hline
 $ \mathcal{H}_{\alpha} $: 
$e_{1}\cdot e_{3}= e_{2},   $ & $ \left(
\begin{array}{ccc}
0 & a_{12} & 0 \\
a_{21} & a_{22} & 0 \\
0 & 0 & 0 \\  
\end{array}
\right) $;$ \left(
\begin{array}{ccc}
0 & 0 & 0 \\
0 & a_{22} & a_{23} \\
0 & a_{32} & 0 \\  
\end{array}
\right)$ $ a_{23}\neq 0 $  \cr
$ e_{3}\cdot e_{1}= \alpha e_{2,}$ & $ \left(
\begin{array}{ccc}
0 & 0 & 0 \\
0 & a_{22} & a_{23} \\
0 & a_{32} & a_{33} \\  
\end{array}
\right) $ $a_{33}\neq 0$    \cr
$\alpha \in  \mathbb{C}-\lbrace 1 \rbrace  $ & $ \left(
\begin{array}{ccc}
0 & a_{12} & 0 \\
a_{21} & a_{22} & 0 \\
0 & a_{32} & 0 \\  
\end{array}
\right) $ $ a_{21}= a_{12}(1 + \alpha),$ $ a_{32}\neq 0 $\cr
&$ \left(
\begin{array}{ccc}
a_{11} & a_{12} & 0 \\
a_{21} & a_{22} & 0 \\
0 & 0 & 0 \\  
\end{array}
\right)$ $ a_{11}\neq 0,$ $\alpha= 0$\cr
& $ \left(
\begin{array}{ccc}
a_{11} & a_{12} & 0 \\
a_{21} & a_{22} & 0 \\
0 & 0 & 0 \\  
\end{array}
\right)$ $  \alpha= -1$ \cr
 & $ \left(
\begin{array}{ccc}
a_{11} & a_{12} & 0 \\
a_{21} & a_{22} & 0 \\
0 & 0 & 0 \\  
\end{array}
\right)$ $a_{11} \neq 0,$ $ \alpha \neq 0, -1  $\cr
&$ \left(
\begin{array}{ccc}
0 & 0 & 0 \\
a_{21} & a_{22} & 0 \\
a_{31} & a_{32} & 0 \\  
\end{array}
\right)$ $ a_{31}\neq 0,$ $\alpha = 0$\cr 
&$ \left(
\begin{array}{ccc}
0 & a_{12} & 0 \\
0 & a_{22} & a_{13} \\
0 & a_{32} & 0 \\  
\end{array}
\right)$  
$ a_{23}, a_{12}\neq 0 $, $ a_{32}= \frac{1}{1 + \alpha} a_{23}, \alpha \neq -1$\cr
& $ \left(
\begin{array}{ccc}
a_{11} & 0 & a_{13} \\
0 & a_{22} & a_{23} \\
a_{31} & 0 & a_{33} \\  
\end{array}
\right)$ $ a_{11}\neq 0 $, $ a_{31}\neq 0 $, $ a_{31}=a_{13}\neq 0, $ \cr
&   $ a_{33}= \frac{a_{31}^{2}}{a_{11}}$ \cr
&$ \left(
\begin{array}{ccc}
a_{11} & 0 & a_{13} \\
a_{21} & a_{22} & a_{23} \\
a_{31} & a_{32} & a_{33} \\  
\end{array}
\right)$ $a_{21}, a_{11}, a_{31}, a_{13}\neq 0 $,$  a_{31}= a_{13},$ \cr
& $ a_{32}=\frac{a_{33}a_{21}}{a_{31}}, a_{33}= \frac{a_{31}^{2}}{a_{11}}$ \cr
&
$ \left(
\begin{array}{ccc}
a_{11} & a_{12} & a_{13} \\
a_{21} & a_{22} & a_{23} \\
a_{31} & a_{32} & a_{33} \\  
\end{array}
\right)$ $ a_{21}, a_{11}, a_{31}, a_{12}, a_{13}\neq 0, $ \cr
& $ a_{32}=\frac{a_{33}(a_{21} - a_{12})}{a_{31}} $, $a_{33}=\frac{a_{31}^{2}}{a_{11}}$\cr
& $\left(
\begin{array}{ccc}
0 & a_{12} & 0 \\
a_{21} & a_{22} & a_{23} \\
0 & a_{32} & a_{33} \\  
\end{array}
\right) $ $ a_{33},  a_{21},  a_{12} \neq 0 $, $  a_{12}=\frac{ a_{21}(1 + \alpha)}{\alpha - 1}, $ \cr
& $  a_{23}=\frac{a_{33}(\alpha a_{21} -  a_{12}) }{ a_{21}},  a_{32}=\frac{ a_{23}( a_{12}- \alpha  a_{21})}{ a_{12}(1+ \alpha)- \alpha  a_{21}}$
\\ \hline
\end{tabular}
\end{center}
\end{proposition}

\textbf{Proof: }
Taking $ a_{11}= a_{22}= a_{33}= 0$ and $ a_{12}= -a_{21},$ $ a_{23}= -a_{32}$ and $ a_{13}= -a_{31} $   the equation $ r_{12}r_{13} + r_{13}r_{23} - r_{23}r_{12}= 0$ provides the following equations: 
\beqs
&&a_{13}= a_{31}=0; \cr
&& a_{12}a_{32} + \alpha a_{32}a_{12} + \alpha a_{23}a_{21}= a_{12}a_{23} + \alpha a_{32}a_{21}.
\eeqs
and the antisymmetric solutions $ r_{1} $ and $ r_{2} $ follow.

$ \hfill \square $ 
\subsection{Antisymmetric infinitesimal Heisenberg bialgebras}
\begin{proposition}
Let us consider the linear map $ \Delta: \mathcal{H}_{\alpha}\rightarrow \mathcal{H}_{\alpha}\otimes \mathcal{H}_{\alpha}  $ defined by 
\beqs
\Delta_{t}(x)= (id\otimes L(x) - R(x)\otimes id)r_{t}, t=1, 2;
\eeqs
for all $ x \in   \mathcal{H}_{\alpha} $. Then, ($  \mathcal{H}_{\alpha}, \Delta_{t} $) is an antisymmetric infinitesimal bialgebra.
\end{proposition}
\textbf{Proof: }
It stems from Corollary \ref{c2}.

$ \hfill \square $

\begin{corollary}
Let $ \mathcal{H}^{\ast}_{\alpha} $  be the dual space of $ \mathcal{H}_{\alpha} $. There exists an associative algebra structure 
$ " \circ " $ on the dual space $ \mathcal{H}^{\ast}_{\alpha} $ given by the dual (linear) map 
$ \Delta^{\ast}_{t} : \mathcal{H}^{\ast}_{\alpha} \otimes \mathcal{H}^{\ast}_{\alpha} \rightarrow \mathcal{H}^{\ast}_{\alpha} $ of $ \Delta_{t} $ defined by
\beq
&&e^{\ast}_{i} \circ e^{\ast}_{j}= - \delta_{ik} \delta_{jk}(\alpha + 1)a_{12} e^{\ast}_{3}, \ \
i, j=1, 2, 3; \ k = 2;\ for \ t= 1;\cr
&&e^{\ast}_{i} \circ e^{\ast}_{j}=  \delta_{ik} \delta_{jk}(\alpha + 1)a_{23} e^{\ast}_{1}, \ \
i, j=1, 2, 3; \ k = 2;\ for \ t= 2.
\eeq
Then, 
 $ (\mathcal{H}_{\alpha}, \mathcal{H}^{\ast}_{\alpha}, R^{\ast}_{\cdot}, L^{\ast}_{\cdot},  R^{\ast}_{\circ}, L^{\ast}_{\circ}) $
 is a matched pair of associative algebras.
\end{corollary}
\textbf{Proof: }
We have
\beqs 
&&\Delta_{1}(e_{1}) = 0 ,  \Delta_{1}(e_{2}) = 0 ,  \Delta_{1}(e_{3}) = -(\alpha + 1)a_{12} e_{2} \otimes e_{2}; \cr
&&\Delta_{2}(e_{1}) =(\alpha + 1)a_{23} e_{2} \otimes e_{2} ,  \Delta_{2}(e_{2}) = 0 ,  \Delta_{2}(e_{3}) =0.
\eeqs
Then, we deduce the following results
\beqs
&&e^{\ast}_{i} \circ e^{\ast}_{j}= - \delta_{ik} \delta_{jk}(\alpha + 1)a_{12} e^{\ast}_{3}, \ \
i, j=1, 2, 3; \ k = 2;\ for \ t= 1;\cr
&&e^{\ast}_{i} \circ e^{\ast}_{j}=  \delta_{ik} \delta_{jk}(\alpha + 1)a_{23} e^{\ast}_{1}, \ \
i, j=1, 2, 3; \ k = 2;\ for \ t= 2.
\eeqs
Hence,  the operation $ \circ $ is associative on $ \mathcal{H}^{\ast}_{\alpha} $. Now, show that
 $ (\mathcal{H}_{\alpha}, \mathcal{H}^{\ast}_{\alpha}, R^{\ast}_{\cdot}, L^{\ast}_{\cdot},  R^{\ast}_{\circ}, L^{\ast}_{\circ}) $ 
is a matched pair of associative algebras. Consider the solution $ r_{1}  $. Set 
\beqs
&& x = x_{1}e_{1} + x_{2}e_{2} + x_{3}e_{3}, \ \ y = y_{1}e_{1} + y_{2}e_{2} + y_{3}e_{3}, \cr
&& a^{\ast} = a_{1}e^{\ast}_{1} + a_{2}e^{\ast}_{2} + a_{3}e^{\ast}_{3}, \ b^{\ast} = b_{1}e^{\ast}_{1} + b_{2}e^{\ast}_{2} + b_{3}e^{\ast}_{3}, 
\eeqs
for all $ x, y \in \mathcal{H}_{\alpha} $ and $ a^{\ast}, b^{\ast} \in \mathcal{H}^{\ast}_{\alpha} $. From $ e_{1}\cdot e_{3}= e_{2};  e_{3}\cdot e_{1}= \alpha e_{2} $ and 
 $ e^{\ast}_{2}\circ e^{\ast}_{2} = -(\alpha + 1)a_{12}e^{\ast}_{3} $, respectively,  we have 
\beqs
R^{\ast}_{\cdot}(e_{3})e^{\ast}_{2} &=& e^{\ast}_{1}, \ \ L^{\ast}_{\cdot}(e_{1})e^{\ast}_{2} = e^{\ast}_{3},  \cr
R^{\ast}_{\cdot}(e_{1})e^{\ast}_{2} &=& \alpha e^{\ast}_{3}, \ \ L^{\ast}_{\cdot}(e_{3})e^{\ast}_{2} = \alpha e^{\ast}_{1},  \cr
R^{\ast}_{\circ}(e^{\ast}_{2})e_{3} &=& -(\alpha + 1)a_{12} e_{2}, \ \ L^{\ast}_{\circ}(e^{\ast}_{2})e_{3} = -(\alpha + 1)a_{12} e_{2}.
\eeqs
Then,  we obtain:
\beqs
&&R^{\ast}_{\cdot}(x)(a^{\ast}\circ b^{\ast})=-a_{2}b_{2}(\alpha + 1)a_{12} R^{\ast}_{\cdot}(x)(e^{\ast}_{3})= 0, \cr
&&R^{\ast}_{\cdot}(L^{\ast}_{\circ}(a^{\ast})x)b^{\ast}=-a_{2}x_{3}(\alpha + 1)a_{12}  R^{\ast}_{\cdot}(e_{2})b^{\ast}= 0, \cr
&&(R^{\ast}_{\cdot}(x)a^{\ast})\circ b^{\ast}= (x_{3}a_{2}e^{\ast}_{1}+ x_{1}a_{2}\alpha e^{\ast}_{3})\circ b^{\ast}= 0.
\eeqs
 Thus,   the equation (\ref{eq16}) is satisfied and
\beqs
&&R^{\ast}_{\cdot}(R^{\ast}_{\circ}(a^{\ast})x)b^{\ast}=-a_{2}x_{3}(\alpha + 1)a_{12} R^{\ast}_{\cdot}(e_{2})b^{\ast}= 0, \cr
&&L^{\ast}_{\cdot}(x)a^{\ast}\circ b^{\ast}= (x_{1}a_{2}e^{\ast}_{3} + x_{3}a_{2}\alpha e^{\ast}_{1})\circ b^{\ast}= 0, \cr
&&L^{\ast}_{\cdot}[L^{\ast}_{\circ}(b^{\ast})x]a^{\ast}=-b_{2}x_{3} (\alpha +1)a_{12}L^{\ast}_{\cdot}(e_{2})a^{\ast}= 0, \cr
&&a^{\ast}\circ(R^{\ast}_{\cdot}(x)b^{\ast})= a^{\ast}\circ (x_{3}b_{2}e^{\ast}_{1} + x_{1}b_{2}\alpha e^{\ast}_{3})= 0 
\eeqs
 implying that
 the equation (\ref{eq17}) holds. Hence,  from the  Theorem \ref{pro2.2.2}, we conclude that \\ 
$ (\mathcal{H}_{\alpha}, \mathcal{H}^{\ast}_{\alpha}, R^{\ast}_{\cdot}, L^{\ast}_{\cdot},  R^{\ast}_{\circ}, L^{\ast}_{\circ}) $ 
is a matched pair of associative algebras.
$ \hfill \square $
\subsection{Heisenberg Frobenius algebras}

\begin{theorem}
$(\mathcal{H}_{\alpha}, \cdot)$ and $(\mathcal{H}^{\ast}_{\alpha}, \circ)$ are two associative algebras. 
There is a double construction of a symmetric Frobenius algebra, denoted by 
$(\mathcal{H}_{\alpha} \bowtie^{R^{\ast}_{\cdot}, L^{\ast}_{\cdot}}_{ R^{\ast}_{\circ}, L^{\ast}_{\circ}} \mathcal{H}^{\ast}_{\alpha}, \mathcal{B})$
 or simply $(\mathcal{H}_{\alpha} \bowtie \mathcal{H}^{\ast}_{\alpha}, \mathcal{B}),$ where $ \mathcal{B} $ is the symmetric bilinear form on the direct sum
 $\mathcal{H}_{\alpha} \oplus \mathcal{H}^{\ast}_{\alpha} $ of the underlying vector spaces of $ \mathcal{H}_{\alpha} $ and $ \mathcal{H}^{\ast}_{\alpha} $ given by (\ref{eq2}). 
The 
 product denoted by $ \ast $  on $ \mathcal{H}_{\alpha} \oplus \mathcal{H}^{\ast}_{\alpha} $ defined by \beqs
(x + a^{\ast})\ast(y + b^{\ast}) &=& (x\cdot y + R^{\ast}_{\circ}(a^{\ast})y + L^{\ast}_{\circ}(b^{\ast})x) + (a^{\ast} \circ b^{\ast} + R^{\ast}_{\cdot}(x)b^{\ast} + L^{\ast}_{\cdot}(y)a^{\ast}) \cr
                                  &&   \mbox{ for all } x, y 
\in  \mathcal{H}, a^{\ast}, b^{\ast} \in  \mathcal{H}^{\ast} 
\eeqs 
 is given in Table 3.

  \begin{center}
{\bf Table 3}: 
 Antisymmetric solutions and Frobenius algebra structures.
\begin{tabular}{c c | c c}
\hline
Associative & Antisymmetric  & Frobenius algebra \cr
 algebra  & solutions     &  structures 
    \\ \hline
$ \mathcal{H}_{\alpha}$  & $ a_{12}e_{1}\otimes e_{2} +  $ & $ e_{1}\ast e_{3}= e_{2}; e_{3}\ast e_{1}= \alpha e_{2};$ $ e^{\ast}_{2}\ast e_{1}= e^{\ast}_{3};$ $ e_{1}\ast e^{\ast}_{2}=\alpha e^{\ast}_{3}$ \cr
& $ a_{21}e_{2}\otimes e_{1}$  & $ e^{\ast}_{2}\ast e^{\ast}_{2}= -(1+ \alpha)a_{12}e^{\ast}_{3};$
   $ e_{3}\ast e^{\ast}_{2}= e^{\ast}_{1} - (\alpha + 1)a_{12}e_{2} $ \cr
  &  $ a_{12} = - a_{21} $ & $ e^{\ast}_{2}\ast e_{3}= \alpha e^{\ast}_{1} -(\alpha + 1)a_{12}e_{2}$    \cr
\cline{2-3}  &  $  a_{23}e_{2}\otimes e_{3} + $  &  $ e_{1}\ast e_{3}= e_{2}$; $ e_{3}\ast e_{1}= \alpha e_{2}$; $ e_{3}\ast e^{\ast}_{2}= e^{\ast}_{1};$ $ e^{\ast}_{2}\ast e_{3}= \alpha e^{\ast}_{1}  $  \cr
 &  $ a_{32}e_{3}\otimes e_{2} $  & $ e^{\ast}_{2}\ast e^{\ast}_{2}= (1 + \alpha)a_{23}e^{\ast}_{1} $; $ e_{1}\ast e^{\ast}_{2}= \alpha e^{\ast}_{3} + (1 + \alpha)a_{23}e_{2} $  \cr
                   & $ a_{32}= -a_{23} $ &  $  e^{\ast}_{2}\ast e_{1}= e^{\ast}_{3} + (1 + \alpha)a_{23}e_{2}  $ 
 \\ \hline
\end{tabular}
\end{center}
\end{theorem}

\begin{remark}
For $ \alpha= 0 $, we obtain the double constructions of  Heisenberg Frobenius  algebras. 
\end{remark}

\begin{theorem}
$(\mathcal{H}_{0}, \Delta_{t})$ is an antisymmetric infinitesimal bialgebra.  
Then, there is a canonical antisymmetric infinitesimal bialgebra structure on the direct sum $ \mathcal{H}_{0} \oplus \mathcal{H}^{\ast}_{0} $ of the underlying vector spaces of $ \mathcal{H}_{0} $ and $ \mathcal{H}^{\ast}_{0} $.
\end{theorem}
\textbf{Proof: }
Take $a_{12}=1$. It stems step by step from the proof of Theorem 2.3.6 in \cite{[C.Bai5]}. In fact, let 
$ r \in \mathcal{H}_{0}\otimes \mathcal{H}^{\ast}_{0} \subset (\mathcal{H}_{0} \oplus \mathcal{H}^{\ast}_{0})\otimes  (\mathcal{H}_{0} \oplus \mathcal{H}^{\ast}_{0})$ 
correspond to the identity map $ id: \mathcal{H}_{0} \oplus \mathcal{H}^{\ast}_{0} \rightarrow \mathcal{H}_{0} \oplus \mathcal{H}^{\ast}_{0}$. Then,
 $ r = e_{1}\otimes e^{\ast}_{1} + e_{2}\otimes e^{\ast}_{2} + e_{3}\otimes e^{\ast}_{3} $. Consider the associative algebra structure 
$ " \ast " $ on  $ \mathcal{H} \oplus \mathcal{H}^{\ast} $ given by $ \mathcal{H}\mathcal{D}(\mathcal{H}) $. Next,  let us show that, if 
$ r $ satisfies the equations $(22)$ and $(23)$ in \cite{[C.Bai5]}, then 
\beqs
\Delta_{ \mathcal{H}\mathcal{D}( \mathcal{H})} = (id \otimes L(u) - R(u)\otimes id)r,
\eeqs
for all $ u \in  \mathcal{H}\mathcal{D}( \mathcal{H}) $,  induces an antisymmetric infinitesimal bialgebra structure on 
$  \mathcal{H}\mathcal{D}( \mathcal{H}) $. We have
\beqs
(L(e_{1})\otimes id - id\otimes R(e_{1}))(r + \sigma(r)) &=& e_{2}\otimes e^{\ast}_{3} - e_{2}\otimes e^{\ast}_{3} = 0, \cr
(L(e_{2})\otimes id - id\otimes R(e_{2}))(r + \sigma(r)) &=& 0, \cr
(L(e_{3})\otimes id - id\otimes R(e_{3}))(r + \sigma(r)) &=& -e_{2}\otimes e_{2} + e^{\ast}_{1}\otimes e_{2} + e_{2}\otimes e_{2} \cr  
                                                  && -e^{\ast}_{1}\otimes e_{2} = 0, \cr
(L(e^{\ast}_{1})\otimes id - id\otimes R(e^{\ast}_{1}))(r + \sigma(r)) &=& 0, \cr
(L(e^{\ast}_{2})\otimes id - id\otimes R(e^{\ast}_{2}))(r + \sigma(r)) &=& e^{\ast}_{3}\otimes e^{\ast}_{1} - e_{2}\otimes e^{\ast}_{3} - 
e^{\ast}_{3}\otimes e_{2} + e_{2}\otimes e^{\ast}_{3}  \cr
&& + e^{\ast}_{3}\otimes e_{2} - e^{\ast}_{3}\otimes e^{\ast}_{1}  \cr
&=& 0, \cr
(L(e^{\ast}_{3})\otimes id - id\otimes R(e^{\ast}_{3}))(r + \sigma(r)) &=& 0. 
\eeqs

Hence, $ r $    well satisfies the equation $(23)$. Further,  we have  the relations
\beqs
r_{12} &=& e_{1}\otimes e^{\ast}_{1}\otimes 1 + e_{2}\otimes e^{\ast}_{2}\otimes 1 + e_{3}\otimes e^{\ast}_{3}\otimes 1, \cr
r_{13} &=& e_{1}\otimes 1\otimes  e^{\ast}_{1} + e_{2}\otimes 1\otimes  e^{\ast}_{2} + e_{3}\otimes 1\otimes  e^{\ast}_{3}, \cr
r_{23} &=& 1 \otimes e_{1}\otimes e^{\ast}_{1} + 1 \otimes e_{2}\otimes e^{\ast}_{2} + 1 \otimes e_{3}\otimes e^{\ast}_{3}, \cr
\eeqs
such that  
\beqs
r_{12}r_{13} &=& e_{2}\otimes e^{\ast}_{1}\otimes e^{\ast}_{3}, \cr
r_{13}r_{23} &=& - e_{2}\otimes e_{2}\otimes e^{\ast}_{3}, \cr
r_{23}r_{12} &=& - e_{2}\otimes e_{2}\otimes e^{\ast}_{3} + e_{2}\otimes e^{\ast}_{1}\otimes e^{\ast}_{3} 
\eeqs
 providing   $ r_{12}r_{13} + r_{13}r_{23} - r_{23}r_{12} = 0 $. Therefore, $ r $ satisfies the equation (22).
$ \hfill \square $

We denote such  a canonical antisymmetric infinitesimal bialgebra  by   $ \mathcal{H}_{0}\mathcal{D}(\mathcal{H}_{0}) ,$ and call it   \textbf{Heisenberg associative double} (with respect to   the name given 
to such a Frobenius bialgebra in \cite{[C.Bai5]} in 2010).

\begin{theorem} Given
$ r_{1}$ and $ r_{2}$   the antisymmetric solutions of the associative Yang-Baxter equation in Heisenberg associative algebra  $ \mathcal{H}_{0} ,$ then, the associative algebra structures on the associative double $ \mathcal{H}_{0}\mathcal{D}( \mathcal{H}_{0}) $ are given 
from the products  in $  \mathcal{H}_{0},$ respectively, as follows:
\beqs 
&&e^{\ast}_{2}\circ e^{\ast}_{2} = L^{\ast}_{\cdot}(r_{1}(e^{\ast}_{2}))e^{\ast}_{2}, \cr
&&e_{3}\ast e^{\ast}_{2} = R^{\ast}_{\cdot}(e_{3})e^{\ast}_{2} - r_{1}(R^{\ast}_{\cdot}(e_{3})e^{\ast}_{2}), \cr
&&e^{\ast}_{2}\ast e_{3} = r_{1}(e^{\ast}_{2})\cdot e_{3},\cr
&&e^{\ast}_{2}\ast e_{1} = L^{\ast}_{\cdot}(e_{1})e^{\ast}_{2}.
\eeqs
\beqs 
&&e^{\ast}_{2}\circ e^{\ast}_{2} = R^{\ast}_{\cdot}(r_{2}(e^{\ast}_{2}))e^{\ast}_{2}, \cr
&&e^{\ast}_{2}\ast e_{1} = L^{\ast}_{\cdot}(e_{1})e^{\ast}_{2} - r_{2}(L^{\ast}_{\cdot}(e_{1})e^{\ast}_{2}), \cr
&&e_{1}\ast e^{\ast}_{2} = e_{1}\cdot r_{2}(e^{\ast}_{2}),\cr
&&e_{3}\ast e^{\ast}_{2} = R^{\ast}_{\cdot}(e_{3})e^{\ast}_{2}.
\eeqs

Therefore,  $ r_{1} $ and $ r_{2} $ can be regarded as of the linear maps from $  \mathcal{H}^{\ast}_{0} $ to $  \mathcal{H}_{0} $ in the following way, respectively:
\beqs
r_{1}(e^{\ast}_{2}) = -a_{12}e_{1}, \ \ r_{1}(e^{\ast}_{1}) = a_{12}e_{2}, \ \ r_{1}(e^{\ast}_{3}) = 0.
\eeqs
\beqs
r_{2}(e^{\ast}_{2}) = a_{23}e_{3}, \ \ r_{2}(e^{\ast}_{3}) = -a_{23}e_{2}, \ \ r_{2}(e^{\ast}_{1}) = 0.
\eeqs

Moreover,  $ r_{1} $ and $ r_{2} $ are compatible with the  associative algebras structures  on $  \mathcal{H}_{0} $ and $  \mathcal{H}^{\ast}_{0} $.
\end{theorem}
%
%

\begin{theorem}
For all $ a^{\ast}, b^{\ast} \in  \mathcal{H}^{\ast}_{0} $, $ r_{t}, t=1, 2 $ satisfies
\beqs 
r_{t}(a^{\ast})\cdot r_{t}(b^{\ast}) = r_{t}(R^{\ast}_{\cdot}(r_{t}(a^{\ast}))b^{\ast} + L^{\ast}_{\cdot}(r_{t}(b^{\ast}))a^{\ast}).
\eeqs 
Therefore,  $ r_{t} $ is an $ \mathcal{O}- $operator associated to the bimodule $(R^{\ast}_{\cdot}, L^{\ast}_{\cdot})$.
\end{theorem}
\textbf{Proof: }
We have:
 \beqs
 r_{1}(a^{\ast})&=& r_{1}(a_{1}e^{\ast}_{1} + a_{2}e^{\ast}_{2} + a_{3}e^{\ast}_{3}) \cr
            &=& a_{1}r_{1}(e^{\ast}_{1}) + a_{2}r(e^{\ast}_{2}) + a_{3}r(e^{\ast}_{3}) \cr
            &=&  a_{12}(a_{1}e_{2} - a_{2}e_{1}).
 \eeqs
 Then
\beqs
r_{1}(a^{\ast})\cdot r_{1}(b^{\ast})=0    
\eeqs
and
\beqs
 r_{1}(R^{\ast}_{\cdot}(r_{1}(a^{\ast}))b^{\ast} + L^{\ast}_{\cdot}(r_{1}(b^{\ast}))a^{\ast})
&=&r_{1}(R^{\ast}_{\cdot}(a_{12}(a_{1}e_{2} - a_{2}e_{1}))b^{\ast} + L^{\ast}_{\cdot}(a_{12}(b_{1}e_{2} - b_{2}e_{1}))a^{\ast})\cr
&=&-a_{12}b_{2}a_{2}r_{1}(e^{\ast}_{3})\cr
&=&0.
\eeqs
Thus, we obtain $ r_{1}(a^{\ast})\cdot r_{1}(b^{\ast}) = r_{1}(R^{\ast}_{\cdot}(r_{1}(a^{\ast}))b^{\ast} + L^{\ast}_{\cdot}(r_{1}(b^{\ast}))a^{\ast}) $. Hence, $ r_{1} $ is well an $ \mathcal{O}- $operator associated to the bimodule $(R^{\ast}_{\cdot}, L^{\ast}_{\cdot})$.
 
$ \hfill \square $

\begin{corollary}
$ r_{t} :  \mathcal{H}^{\ast}_{0} \rightarrow  \mathcal{H}_{0} $ is a homomorphism of associative algebras, that is, we have
\beq
r_{t}(a^{\ast} \circ b^{\ast}) = r_{t}(a^{\ast})\cdot r_{t}(b^{\ast}) \mbox{ for all } a^{\ast}, b^{\ast} \in H^{\ast}_{0}.
\eeq
\end{corollary}
\textbf{Proof:}
By direct computation, we obtain the result.
$ \hfill \square $

\begin{theorem}
There exists a dendriform algebra structure on $  \mathcal{H}^{\ast}_{0} $ given by
\beq
a^{\ast} \succ b^{\ast} = R^{\ast}_{\cdot}(r_{t}(a^{\ast}))b^{\ast}, \ \ a^{\ast} \prec b^{\ast} = L^{\ast}_{\cdot}(r_{t}(b^{\ast}))a^{\ast},
\eeq
for all $ a^{\ast}, b^{\ast} \in  \mathcal{H}^{\ast}_{0} $.Then, 
there is an associated associative algebra structure on $  \mathcal{H}^{\ast} $ given by 
\beq
a^{\ast} \ast b^{\ast} = a^{\ast} \prec b^{\ast} + a^{\ast} \succ b^{\ast}.
\eeq  
\end{theorem}
\textbf{Proof:}
It is straightforward.

$ \hfill \square $

\begin{corollary}
$( \mathcal{H}_{0}, \cdot)$ is an associative algebra. Then,
\beq
r = e_{1}\otimes e^{\ast}_{1} -  e^{\ast}_{1}\otimes e_{1} + e_{2}\otimes e^{\ast}_{2} -  e^{\ast}_{2}\otimes e_{2} + e_{3}\otimes e^{\ast}_{3}
 -  e^{\ast}_{3}\otimes e_{3} 
\eeq
is a solution of AYBE
in $  \mathcal{H} \ltimes_{R^{\ast}_{\cdot}, 0}  \mathcal{H}^{\ast} $ or
$  \mathcal{H} \ltimes_{0 ,L^{\ast}}  \mathcal{H}^{\ast} $. Moreover, there is a natural Connes cocycle $ \omega $ on 
$  \mathcal{H} \ltimes_{R^{\ast}_{\cdot}, 0}  \mathcal{H}^{\ast} $ or $  \mathcal{H} \ltimes_{0 ,L^{\ast}}  \mathcal{H}^{\ast} $
 which is given by the equation (\ref{eq7}). 
\end{corollary}
\textbf{Proof: }
We set $ a_{12}=1$. $ (R^{\ast}_{\cdot}, 0,  \mathcal{H}^{\ast}_{0}) $ is a $  \mathcal{H}_{0}- $bimodule, 
then ($\mathcal{H}_{0} \oplus  \mathcal{H}^{\ast}_{0}, \ast $) is an associative algebra 
where the product is defined by
\beqs
(x + a^{\ast})\ast(y + b^{\ast}) = x\cdot y + R^{\ast}_{\cdot}(x)b^{\ast}, \mbox{ for all } x, y 
\in  \mathcal{H}, a^{\ast}, b^{\ast} \in  \mathcal{H}^{\ast}.  
\eeqs 
We denote it by $  \mathcal{H}_{0} \ltimes_{R^{\ast}_{\cdot}, 0}  \mathcal{H}^{\ast}_{0} $ and we simply obtain 
\beqs
       e_{1}\ast e_{3} &=& e_{3}, \cr
e_{3}\ast e^{\ast}_{2} &=& e^{\ast}_{1}, \cr
\eeqs
and $0$  elsewhere. Now, show that   $ r_{12}r_{13} + r_{13}r_{23} - r_{23}r_{12} = 0 $. We have
\beqs
r_{12} &=& e_{1}\otimes e^{\ast}_{1}\otimes 1 -  e^{\ast}_{1}\otimes e_{1}\otimes 1 + e_{2}\otimes e^{\ast}_{2}\otimes 1 -  e^{\ast}_{2}\otimes e_{2}\otimes 1 \cr
&& + e_{3}\otimes e^{\ast}_{3}\otimes 1 -  e^{\ast}_{3}\otimes e_{3}\otimes 1, \cr
r_{13} &=& e_{1}\otimes 1\otimes e^{\ast}_{1} -  e^{\ast}_{1}\otimes 1\otimes e_{1} + e_{2}\otimes 1\otimes e^{\ast}_{2} -  e^{\ast}_{2}\otimes 1\otimes e_{2} \cr
&& + e_{3}\otimes 1\otimes e^{\ast}_{3} -  e^{\ast}_{3}\otimes 1\otimes e_{3},  \cr
r_{23} &=& 1\otimes e_{1}\otimes e^{\ast}_{1} - 1\otimes e^{\ast}_{1}\otimes e_{1} + 1\otimes e_{2}\otimes e^{\ast}_{2} - 1\otimes e^{\ast}_{2}\otimes e_{2} \cr
&& + 1\otimes e_{3}\otimes e^{\ast}_{3} - 1\otimes e^{\ast}_{3}\otimes e_{3}.
\eeqs
Then, 
\beqs
r_{12}r_{13} &=& e_{2}\otimes e^{\ast}_{1}\otimes e^{\ast}_{3} - e^{\ast}_{1}\otimes e^{\ast}_{3}\otimes e_{2}, \cr
r_{13}r_{23} &=& e^{\ast}_{1}\otimes e^{\ast}_{3}\otimes e_{2} - e^{\ast}_{3}\otimes e_{2}\otimes e^{\ast}_{1}, \cr
r_{23}r_{12} &=& - e^{\ast}_{3}\otimes e_{2}\otimes e^{\ast}_{1} + e^{\ast}_{2}\otimes e^{\ast}_{1}\otimes e^{\ast}_{3}, \cr
\eeqs
and  $ r_{12}r_{13} + r_{13}r_{23} - r_{23}r_{12} = 0 $. 
Moreover, show that the antisymmetric bilinear form, denoted by $ \omega $ on the associative algebra 
$   \mathcal{H} \ltimes_{R^{\ast}_{\cdot}, 0}  \mathcal{H}^{\ast},  $ defined by
\beqs
\omega(x + a^{\ast}, y + b^{\ast}) = -\langle x, b^{\ast} \rangle + \langle a^{\ast}, y \rangle \mbox { for all } x, y \in  \mathcal{H}, 
a^{\ast}, b^{\ast} \in  \mathcal{H}^{\ast} ,
\eeqs 
is a Connes cocycle, that is, for all $ z \in  \mathcal{H}, c^{\ast} \in  \mathcal{H}^{\ast}$ 
\beq \label{eq90}
&& \omega[(x + a^{\ast})\ast(y + b^{\ast}),(z + c^{\ast})] + \omega[(y + b^{\ast})\ast(z + c^{\ast}),(x + a^{\ast})] \cr
&& + \omega[(z + c^{\ast})\ast(x + a^{\ast}),(y + b^{\ast})] =0.
\eeq
We have
\beqs
 \omega[(x + a^{\ast})\ast(y + b^{\ast}),(z + c^{\ast})] &=& \omega(x\cdot y + R^{\ast}_{\cdot}(x)b^{\ast}, z + c^{\ast}) \cr
                                                         &=& -\langle x\cdot y, c^{\ast} \rangle + \langle R^{\ast}_{\cdot}(x)b^{\ast}, z \rangle \cr
                                                         &=& - \langle x\cdot y, c^{\ast} \rangle + \langle z\cdot x, b^{\ast} \rangle , \cr \cr
\omega[(y + b^{\ast})\ast(z + c^{\ast}),(x + a^{\ast})] &=& \omega(y\cdot z + R^{\ast}_{\cdot}(y)c^{\ast}, x + a^{\ast}) \cr
                                                         &=& -\langle y\cdot z, a^{\ast} \rangle + \langle R^{\ast}_{\cdot}(y)c^{\ast}, x \rangle \cr
                                                         &=& - \langle y\cdot z, a^{\ast} \rangle + \langle x\cdot y, c^{\ast} \rangle , \cr \cr
\omega[(z + c^{\ast})\ast(x + a^{\ast}),(y + b^{\ast})] &=& \omega(z\cdot x + R^{\ast}_{\cdot}(z)a^{\ast}, y + b^{\ast}) \cr
                                                         &=& -\langle z\cdot x, b^{\ast} \rangle + \langle R^{\ast}_{\cdot}(z)a^{\ast}, y \rangle \cr
                                                         &=& - \langle z\cdot x, b^{\ast} \rangle + \langle y\cdot z, a^{\ast} \rangle ,                                                                          
\eeqs
such that  the equation (\ref{eq90}) is satisfied.

$ \hfill \square $

\section{Double constructions of  Heisenberg Connes cocycles}
In this section, due to the complexity of the computations, we only determine  important  compatible dendriform algebras on the extended associative Heisenberg algebra $\mathcal{H}_{\alpha}$. Besides, we perform the classification of solutions of the  $D-$equations in them, 
and establish the corresponding double constructions of Connes cocycles.

We start by considering  the following simpler case.
\begin{proposition}
Let $ \prec $  and $ \succ $  be two operations on $\mathcal{H}_{0}$. The  structure ($\mathcal{H}_{0}, \succ, \prec$)  defined by :
\beq \label{D1}
&&e_{1}\prec e_{3}= (1-\lambda)e_{2}, \cr
&&e_{1}\succ e_{3}=\lambda e_{2}
\eeq
is a  dendriform algebra structure.
\end{proposition}
\textbf{Proof} 
It follows from a direct computation.

$ \hfill \square $
\subsubsection{  Classification of the solutions of $ D- $equation in ($\mathcal{H}_{0}, \succ, \prec$) }
\begin{theorem}\label{t3}
Let $ r = \sum_{i, j}a_{ij}e_{i}\otimes e_{j} \in \mathcal{H}_{0}\otimes \mathcal{H}_{0} $, $ i, j=1, 2, 3 $. We have following solutions: 
\benum
\item[(i)] 
\beqs
\lbrace
\left(
\begin{array}{ccc}
a_{11} & a_{12} & 0 \\
a_{21} & a_{22} & 0 \\
0 & 0 & 0 \\  
\end{array}
\right), a_{11} \neq 0, \lambda= 0
\rbrace ;
\eeqs
\item[(ii)]
 \beqs
\lbrace
\left(
\begin{array}{ccc}
0 & a_{12} & a_{13} \\
0 & a_{22} & a_{23} \\
0 & 0 & 0 \\  
\end{array}
\right), a_{13} \neq 0, \lambda= 0
			\rbrace ;
			\eeqs
\item[(iii)]
\beqs
\lbrace
\left(
\begin{array}{ccc}
0 & a_{12} & 0 \\
a_{21} & a_{22} & a_{23} \\
0 & a_{32} & 0 \\  
\end{array}
\right), a_{32} \neq 0,  \lambda= 0
\rbrace ;
\eeqs
\item[(iv)]
\beqs
\lbrace
\left(
\begin{array}{ccc}
0 & a_{12} & 0 \\
0 & a_{22} & a_{23} \\
0 & 0 & 0 \\  
\end{array}
\right), \lambda= 0
\rbrace ;
\eeqs
\item[(v)]
\beqs
\lbrace
\left(
\begin{array}{ccc}
0 & a_{12} & 0 \\
a_{21} & a_{22} & 0 \\
0 & 0 & 0 \\  
\end{array}
\right), a_{21}\neq 0, \lambda= 0
\rbrace ;
\eeqs
\item[(vi)]
\beqs
\lbrace
\left(
\begin{array}{ccc}
0 & 0 & 0 \\
0 & a_{22} & a_{23} \\
0 & a_{32} & a_{33} \\  
\end{array}
\right), a_{33} \neq 0, \lambda= 0
\rbrace ;
\eeqs
\item[(vii)]
\beqs
\lbrace
\left(
\begin{array}{ccc}
0 & 0 & 0 \\
a_{21} & a_{22} & 0 \\
a_{31} & a_{32} & 0 \\  
\end{array}
\right), a_{31}\neq 0,  \lambda= 0
\rbrace ;
\eeqs
\item[(viii)]
\beqs
\lbrace
\left(
\begin{array}{ccc}
a_{11} & a_{12} & 0 \\
a_{21} & a_{22} & 0 \\
0 & 0 & 0 \\  
\end{array}
\right), a_{11}\neq 0, \lambda \neq 0
\rbrace ;
\eeqs
\item[(ix)]
\beqs
\lbrace
\left(
\begin{array}{ccc}
0 & a_{12} & 0 \\
a_{21} & a_{22} & a_{23} \\
0 & a_{32} & 0 \\  
\end{array}
\right), a_{12} \neq 0, \lambda \neq 0, a_{32}=\frac{a_{21}a_{23}(1- \lambda)}{a_{12}} + \lambda a_{23}
\rbrace ;
\eeqs
\item[(x)]
\beqs
\lbrace
\left(
\begin{array}{ccc}
0 & 0 & 0 \\
a_{21} & a_{22} & a_{23} \\
a_{31} & a_{32} & a_{33} \\  
\end{array}
\right), \lambda\neq 0, \lambda= 1
\rbrace ;
\eeqs
\item[(xi)]
\beqs
\lbrace
\left(
\begin{array}{ccc}
0 & 0 & 0 \\
a_{21} & a_{22} & 0 \\
a_{31} & a_{32} & 0 \\  
\end{array}
\right), a_{31} \neq 0, \lambda \neq 0, 1
\rbrace ;
\eeqs
\item[(xii)]
\beqs
\lbrace
\left(
\begin{array}{ccc}
0 & 0 & 0 \\
0 & a_{22} & a_{23} \\
0 & a_{32} & a_{33} \\  
\end{array}
\right), \lambda \neq 0, \lambda \neq 1
\rbrace ;
\eeqs
\item[(xiii)]
\beqs
\lbrace
\left(
\begin{array}{ccc}
0 & 0 & 0 \\
a_{21} & a_{22} & 0 \\
0 & a_{32} & 0 \\  
\end{array}
\right), a_{21} \neq 0, \lambda \neq 0, 1
\rbrace .
\eeqs
\eenum
\end{theorem}
\textbf{Proof}
We set $ r = \sum_{i, j}a_{ij}e_{i}\otimes e_{j} \in \mathcal{H}_{0}\otimes \mathcal{H}_{0} $, $ i, j=1, 2, 3 $. Resolving the D-equation $ r_{12}\prec r_{13} + r_{12}\succ r_{13} = r_{12}\prec r_{23} + r_{23}\succ r_{12} $, we obtain the following equations:
\begin{center}
\begin{tabular}{c|c}
  \hline
 $D_{1}:\lambda a_{13}=0$ & $D_{2}: a_{11}a_{31}= 0$ 
\\ \hline
$D_{3}:a_{11}a_{32} -a_{21}a_{13}(1-\lambda)= 0$ & $D_{4}:a_{11}a_{33}= 0$
\\ \hline
$D_{5}: a_{12}a_{31} - \lambda a_{11}a_{23}= 0$ & $D_{6}: a_{12}a_{33} - \lambda a_{13}a_{23}= 0$
\\ \hline
$D_{7}: a_{13}a_{31}= 0$ & $D_{8}: a_{13}a_{33}= 0$
\\ \hline
$D_{9}: a_{11}a_{13}(1-\lambda)= 0$ & $D_{10}: 
a_{12}a_{32} - a_{21}a_{23}(1-\lambda) - \lambda a_{12}a_{23}= 0$
\\ \hline
$D_{11}: a_{13}a_{32} - a_{21}a_{33}(1-\lambda)= 0$ & $D_{12}: a_{11}a_{23}(1-\lambda)= 0$
\\ \hline
$D_{13}: \lambda a_{11}a_{13}= 0$ & $D_{14}: \lambda a_{12}a_{13}= 0$
\\ \hline
$D_{15}: a_{31}a_{33}(1-\lambda)= 0$ & $D_{16}: a_{31}a_{23}(1-\lambda) + \lambda a_{21}a_{33}= 0$
\\ \hline
\end{tabular}
\end{center}
After discussion, we obtain the above solutions.
 
$ \hfill \square $
\subsubsection{Symmetric solutions and Heisenberg Connes cocycle structures}
\begin{proposition}
The double constructions of the Connes cocycle associated with the Heisenberg associative algebra $ \mathcal{H}_{0}$ and his dual space $  \mathcal{H}^{\ast}_{0}$, are given in Table 4.
\begin{center}
{\bf Table 4}: 
 Symmetric solutions and Heisenberg Connes cocycle structures.
\begin{tabular}{c | c}
\hline
Symmetric solutions & Heisenberg Connes cocycle structures   \\ \hline
$ r= a_{11}e_{1}\otimes e_{1} + a_{12}e_{1}\otimes e_{2}  $ & $ e_{1}\ast e_{3}=  e_{2};e^{\ast}_{2}\ast e^{\ast}_{2}= -a_{12}\lambda e^{\ast}_{3} $  \cr
$ + a_{21}e_{2}\otimes e_{1}+ a_{22}e_{2}\otimes e_{2} $; & $ e^{\ast}_{2}\ast e_{1} =\lambda e^{\ast}_{3}; e^{\ast}_{2}\ast e_{3}=- a_{12}e_{2} $  \cr
    $ a_{21}= a_{12} $   & $ e^{\ast}_{2}\ast e^{\ast}_{1} = - a_{11}\lambda e^{\ast}_{3}; e^{\ast}_{1}\ast e_{3} = - a_{11}\lambda e_{2} $ \cr
    & $ e_{3}\ast e^{\ast}_{2} = a_{11}(1-\lambda)e_{1} +(1-\lambda)e^{\ast}_{1}+ (1-\lambda)a_{12}e_{2}   $
\\ \hline
$ r= a_{22}e_{2}\otimes e_{2} + a_{23}e_{2}\otimes e_{3} $ & $ e_{1}\ast e_{3}=  e_{2};e^{\ast}_{2}\ast e^{\ast}_{2}=(\lambda-1)a_{23}e^{\ast}_{1} $  \cr
$ a_{32}e_{3}\otimes e_{2} + a_{33}e_{3}\otimes e_{3} $ & $ e_{1}\ast e^{\ast}_{2}= -a_{23} e_{2}; e_{3}\ast e^{\ast}_{2}=(1-\lambda)e^{\ast}_{1} $  \cr
 $ a_{32}= a_{23} $ & $ e^{\ast}_{3}\ast e^{\ast}_{2} =(\lambda-1)a_{33} e^{\ast}_{1}; e_{1}\ast e^{\ast}_{3}= -a_{33}e_{2} $\cr
 & $  e^{\ast}_{2}\ast e_{1}= \lambda(a_{23}e_{2}+ a_{33}e_{3}+ e^{\ast}_{3}) $ 
  \\ \hline
$ r= a_{12}e_{1}\otimes e_{2} + a_{21}e_{2}\otimes e_{1}  $ & $ e_{1}\ast e_{3}=  e_{2}; e^{\ast}_{2}\ast e^{\ast}_{2}= a_{23}(\lambda-1)e^{\ast}_{1} $  \cr
$ + a_{22}e_{2}\otimes e_{2}+ a_{23}e_{2}\otimes e_{3} $ & $ e^{\ast}_{2}\ast e_{3} = a_{12}(\lambda-1)e_{2}; e_{1}\ast e^{\ast}_{2}= -a_{23}e_{2} $  \cr
    $ + a_{32}e_{3}\otimes e_{2}; a_{21}= a_{12},  $  & $ e_{3}\ast e^{\ast}_{2} =  (1-\lambda)(a_{12}e_{2} + e^{\ast}_{1}) $ \cr
 $ a_{23}= a_{32}$   & $ e^{\ast}_{2}\ast e_{1}= \lambda(a_{23}e_{2} + e^{\ast}_{3} ).$
\\ \hline
\end{tabular}
\end{center}
\end{proposition}
\textbf{Proof:}
Consider the tensor element $ r= a_{12}e_{1}\otimes e_{2} + a_{21}e_{2}\otimes e_{1} + a_{22}e_{2}\otimes e_{2}+ a_{23}e_{2}\otimes e_{3} + a_{32}e_{3}\otimes e_{2}; a_{21}= a_{12}, a_{23}= a_{32}$. r is symmetric and satisfies the $ D- $equation using the equations $ D_{1}- D_{16} $ obtained in the above proof. Using  Corollary \ref{cor1}, we obtain 
\beqs
&&\Delta_{\succ}(e_{1})= -a_{23}e_{2}\otimes e_{2}; \Delta_{\succ}(e_{2})= 0; \Delta_{\succ}(e_{3})= a_{12}(1- \lambda)e_{2}\otimes e_{2};\cr
&&\Delta_{\prec}(e_{1})= a_{23}\lambda e_{2}\otimes e_{2}; \Delta_{\prec}(e_{2})= 0; \Delta_{\prec}(e_{3})= -a_{12}(1- \lambda)e_{2}\otimes e_{2}.
\eeqs
Hence, we obtain on  $ \mathcal{H}^{\ast}_{0} $ the dendriform algebra structure defined by
\beq \label{D2}
e^{\ast}_{2}\succ e^{\ast}_{2}= - a_{23}e^{\ast}_{1} +  a_{12}(1- \lambda)e^{\ast}_{3}; \ e^{\ast}_{2}\prec e^{\ast}_{2}= \lambda a_{23}e^{\ast}_{1} -a_{12}(1- \lambda)e^{\ast}_{3}.
\eeq
Then, from the equations \ref{D1}, \ref{D2}, we have, respectively, the following relations: 
\beqs
&&R^{\ast}_{\prec}(e_{3})e^{\ast}_{2}= (1- \lambda)e^{\ast}_{1};\ L^{\ast}_{\succ}(e_{1})e^{\ast}_{2}= \lambda e^{\ast}_{3};\ L^{\ast}_{\succ}(e^{\ast}_{2})e_{1}= -a_{23}e_{2}; \cr
&&L^{\ast}_{\succ}(e^{\ast}_{2})e_{3}= a_{12}(1- \lambda)e_{2};\ R^{\ast}_{\prec}(e^{\ast}_{2})e_{1}= \lambda a_{23}e_{2}; R^{\ast}_{\prec}(e^{\ast}_{2})e_{3}= -a_{12}(1- \lambda)e_{2}.
\eeqs
Using the associative product on $\mathcal{H}_{0}\oplus \mathcal{H}^{\ast}_{0}$ defined by  
\beqs
(x + a^{\ast})\ast_{\mathcal{H}_{0} \oplus \mathcal{H}_{0}^{\ast}} (y + b^{\ast}) &=&(x \ast_{\mathcal{H}_{0}} y + R^{\ast}_{\prec_{\mathcal{H}_{0}^{\ast}}}(a^{\ast})y 
+ L^{\ast}_{\succ_{\mathcal{H}_{0}^{\ast}}}(b^{\ast})x)\cr 
&&+ (a^{\ast} \ast_{\mathcal{H}_{0}^{\ast}} b^{\ast} + R^{\ast}_{\prec_{\mathcal{H}_{0}}}(x)b^{\ast} + L^{\ast}_{\succ_{\mathcal{H}_{0}}}(y)a^{\ast}), 
\eeqs
we get the results
\beqs
&& e_{1}\ast e_{3}=  e_{2}; e^{\ast}_{2}\ast e^{\ast}_{2}= a_{23}(\lambda-1)e^{\ast}_{1};  
  e^{\ast}_{2}\ast e_{3} = a_{12}(\lambda-1)e_{2}; e_{1}\ast e^{\ast}_{2}= -a_{23}e_{2}   \cr
&& e_{3}\ast e^{\ast}_{2} =  (1-\lambda)(a_{12}e_{2} + e^{\ast}_{1});
  e^{\ast}_{2}\ast e_{1}= \lambda(a_{23}e_{2} + e^{\ast}_{3} ).
\eeqs
$ \hfill \square $ 
\subsection{Some compatible dendriform algebras and classification of solutions of the $D-$equations}
Let $  e_{i}, e_{j}, e_{k} \in \mathcal{H}_{\alpha}; i, j, k=1, 2, 3 $. We have
\beq \label{D2}
&&e_{i}\cdot e_{j}= e_{i}\prec e_{j} + e_{i}\succ e_{j}, \mbox { where we set } \cr
&& e_{i}\succ e_{j} = \sum a^{k}_{ij}e_{k}; e_{i}\prec e_{j} = \sum b^{k}_{ij}e_{k};\cr
&&(e_{i} \prec e_{j}) \prec e_{k} = e_{i} \prec (e_{j} \cdot e_{k}); \cr
&&(e_{i} \succ e_{j}) \prec e_{k} = e_{i} \succ (e_{j} \prec e_{k}); \cr
&& e_{i} \succ (e_{j} \succ e_{k})= ( e_{i} \cdot e_{j}) \succ e_{k}.
\eeq 
Computing the equations (\ref{D2}) with the conditions 
\beq
&&e_{i} \cdot e_{j} = \delta_{ik} \delta_{jl} e_{2}, \ \ i, j=1, 2, 3; \ \ k=1, l= 3;\cr
&&e_{m} \cdot e_{n} =\alpha \delta_{mu} \delta_{nv} e_{2}, \ \ m, n=1, 2, 3; \ \ v=1, u= 3, 
\eeq
 we get the results of Table $5$.
 \subsection{Symmetric solutions and Heisenberg Connes cocycles structures}
 Using the Proposition \ref{pro4.4.6}, we get the results of Table $6$.
\section{Solutions of the associative Yang-Baxter equation on $3$-dimensional non decomposable associative algebras}
 Let $(\mathcal{A}, \cdot)$ be an associative algebra with a basis $\lbrace e_{1}, e_{2}, e_{3} \rbrace$ and $ r = \sum^{3}_{i, j=1} a_{ij}e_{i}\otimes e_{j} \in \mathcal{A} \otimes \mathcal{A}.$ Then, the AYBE (\ref{AYEB}), i.e. $r_{12}r_{13} + r_{13}r_{23} - r_{23}r_{12} = 0,$ is satisfied for
\beqs
r_{12}= \sum^{2}_{i, j=1} a_{ij}e_{i}\otimes e_{j}\otimes 1,
r_{13}= \sum^{2}_{i, j=1} a_{ij}e_{i}\otimes 1\otimes e_{j},
r_{23}= \sum^{2}_{i, j=1} a_{ij} 1\otimes e_{i}\otimes e_{j}.
\eeqs

In the sequel,  we consider the classification of 3-dimensional non decomposable complex associative algebras  from Theorem 4.3 in \cite{[R1]}. Then by direct calculations, we have the following results in the Table $7$.

  \subsection{Antisymmetric solutions and Frobenius algebras structures}
  Using the Proposition \ref{pro2.4.4}, we get the results of Table 8.

{\bf Table 5}:
 Dendriform algebra structures and solutions of $D$-equation.
\begin{center}
\begin{tabular}{c c c}
\hline
 Dendriform  & Solutions of
 \cr
 algebra     &   $D$-equation \cr
 $ \mathcal{A} $    &
    \\ \hline
 $ D^{1}_{\mu_{1}, \mu_{2}, \mu_{3}} $: $e_{1}\prec e_{3}= \mu_{1}e_{2}$ & $ \left(
\begin{array}{ccc}
0 & 0 & 0 \\
a_{21} & a_{22} & 0 \\
0 & 0 & 0 \\  
\end{array}
\right)
$; $ \left(
\begin{array}{ccc}
0 & 0 & 0 \\
0 & a_{22} & 0 \\
0 & a_{32} & 0 \\  
\end{array}
\right)
a_{32}\neq 0    $  \cr
  $  e_{3}\prec e_{1}=\mu_{2}e_{2} $  & $ \left(
\begin{array}{ccc}
0 & 0 & 0 \\
a_{21} & a_{22} & 0 \\
0 & a_{32} & 0 \\  
\end{array}
\right)
a_{32}, a_{21}\neq 0, \alpha= \mu_{2}$ \cr
 $  e_{3}\prec e_{3}=\mu_{3}e_{2} $  & $ \left(
\begin{array}{ccc}
0 & 0 & 0 \\
a_{21} & a_{22} & a_{23} \\
0 & 0 & 0 \\  
\end{array}
\right)a_{23}\neq 0, \mu_{3}=\dfrac{-a_{21}(\mu_{1}+ \mu_{2})}{a_{23}}
$ \cr
$  e_{1}\succ e_{3}=(1-\mu_{1})e_{2} $  & $ \left(
\begin{array}{ccc}
0 & 0 & 0 \\
0 & a_{22} & a_{23} \\
0 & a_{32} & 0 \\  
\end{array}
\right)a_{23}, a_{32}\neq 0, \mu_{3}=\dfrac{-a_{21}(\mu_{1}+ \mu_{2}) + a_{32}\mu_{1}}{a_{23}}
$\cr
$  e_{3}\succ e_{1}=(\alpha- \mu_{2})e_{2}$&$ \left(
\begin{array}{ccc}
0 & 0 & 0 \\
a_{21} & a_{22} & a_{23} \\
0 & a_{32} & 0 \\  
\end{array}
\right)a_{21}\neq 0,
$ \cr
$  e_{3}\succ e_{3}= -\mu_{3}e_{2} $& $\alpha=\dfrac{\mu_{1}a_{23}(a_{32}-a_{21})+ \mu_{2}a_{21}(a_{32}- a_{23}-a^{2}_{23}\mu_{3})}{a_{32}a_{21}}$ \cr
& $ \left(
\begin{array}{ccc}
0 & a_{12} & 0 \\
a_{21} & a_{22} & 0 \\
0 & 0 & 0 \\  
\end{array}
\right)a_{12}\neq 0
$  \cr
&$ \left(
\begin{array}{ccc}
0 & a_{12} & 0 \\
a_{21} & a_{22} & a_{23} \\
0 & 0 & 0 \\  
\end{array}
\right)a_{23}\neq 0, \mu_{3}=\dfrac{-a_{21}(\mu_{1}+ \mu_{2})-a_{12}(1- \mu_{1})}{a_{23}}
$\cr 
&$ \left(
\begin{array}{ccc}
0 & a_{12} & 0 \\
a_{21} & a_{22} & a_{23} \\
0 & a_{32} & 0 \\  
\end{array}
\right)
a_{12}, a_{32}\neq 0, a_{21}\neq a_{12},$ \cr
&$ \alpha=\dfrac{a_{21}a_{23}(\mu_{1} + \mu_{2})+ a^{2}_{23}\mu_{3}-a_{12}a_{23}(1-\mu_{1})-a_{32}(a_{21}\mu_{2}+a_{23}\mu_{1})}{a_{32}(a_{12}- a_{21})}$ \cr
&$ \left(
\begin{array}{ccc}
0 & a_{12} & 0 \\
a_{21} & a_{22} & 0 \\
0 & a_{32} & 0 \\  
\end{array}
\right)a_{12}, a_{32} \neq0, a_{21}=a_{12}, \mu_{2}=-1
$ 
\\ \hline
\end{tabular}
\end{center}

\begin{center}
\begin{tabular}{c c c}
\hline
 Dendriform  & Solutions of
 \cr
 algebra     &   $D$-equation \cr
 $ \mathcal{A} $    &
    \\ \hline
   $ D^{1}_{\mu_{1}, \mu_{2}, \mu_{3}} $:  &$ \left(
\begin{array}{ccc}
0 & a_{12} & 0 \\
a_{21} & a_{22} & a_{23} \\
0 & a_{32} & 0 \\  
\end{array}
\right)a_{12}, a_{32}, a_{23}\neq 0, a_{21}= a_{12},
$\cr
&$\mu_{3}=\dfrac{a_{12}a_{32}(1+ \mu_{2})- a_{12}a_{23}\mu_{2}-a_{22}a_{23}+ a_{32}a_{23}\mu_{1}}{a^{2}_{23}}  $\cr 
&$ \left(
\begin{array}{ccc}
0 & 0 & 0 \\
a_{21} & a_{22} & 0 \\
a_{31} & a_{32} & 0 \\  
\end{array}
\right)a_{31}\neq 0
$ \cr
&$ \left(
\begin{array}{ccc}
0 & 0 & 0 \\
a_{21} & a_{22} & a_{23} \\
a_{31} & a_{32} & 0 \\  
\end{array}
\right)a_{23}\neq 0, \mu_{1}=\mu_{2}=\mu_{3}=0
$ \cr
&$ \left(
\begin{array}{ccc}
0 & a_{12} & 0 \\
a_{21} & a_{22} & a_{23} \\
a_{31} & a_{32} & 0 \\  
\end{array}
\right)a_{31}, a_{12}, a_{23}, \mu_{3}\neq 0, 
$ \cr
&$a_{12}=- \mu_{3}a_{23}, a_{32}=\dfrac{a_{23}(a_{23}\mu_{3}+ a_{12}-a_{32}\mu_{3})}{a_{12}} $\cr
&$ \left(
\begin{array}{ccc}
0 & 0 & 0 \\
0 & a_{22} & a_{23} \\
0 & a_{32} & a_{33} \\  
\end{array}
\right)a_{33}\neq 0, \mu_{3}=0
$\cr 
  & $ \left(
\begin{array}{ccc}
0 & 0 & 0 \\
a_{21} & a_{22} & a_{23} \\
0 & 0 & a_{33} \\  
\end{array}
\right)a_{33}, a_{21}\neq 0, \mu_{1}, \mu_{3}=0
$\cr
&$ \left(
\begin{array}{ccc}
0 & 0 & 0 \\
a_{21} & a_{22} & a_{23} \\
0 & a_{32} & a_{33} \\  
\end{array}
\right)a_{33}, a_{21}, a_{32}\neq 0, \mu_{1} ,\mu_{3}=0, \mu_{2}=\alpha
$ \cr
&$ \left(
\begin{array}{ccc}
0 & a_{12} & 0 \\
a_{21} & a_{22} & a_{23} \\
0 & a_{32} & a_{33} \\  
\end{array}
\right)a_{33}, a_{12}, a_{23}, a_{21}\neq 0, \mu_{3}=0,
$ \cr
&$\alpha= \dfrac{a_{21}\mu_{1}}{a_{12}}, \mu_{1}=\dfrac{a_{12}a_{32}(1+\alpha)-a_{32}a_{21}(\mu_{2}-\alpha)}{a_{21}a_{23}} $\cr
&$ \left(
\begin{array}{ccc}
0 & a_{12} & 0 \\
a_{21} & a_{22} & 0 \\
0 & 0 & a_{33} \\  
\end{array}
\right)a_{33}, a_{21}, a_{12}\neq 0, \alpha=\dfrac{a_{21}\mu_{1}}{a_{12}}, \mu_{3}=0
$ \cr
&$ \left(
\begin{array}{ccc}
0 & a_{12} & 0 \\
a_{21} & a_{22} & 0 \\
0 & a_{32} & a_{33} \\  
\end{array}
\right)a_{33}, a_{32}, a_{21}, a_{12}\neq 0, a_{12}\neq -a_{21}, \mu_{3}=0, $ \cr
&$\alpha= \dfrac{a_{21}\mu_{1}}{a_{12}}=\dfrac{a_{21}\mu_{2}- a_{12}}{a_{12}+a_{21}}  $ \cr
&$ \left(
\begin{array}{ccc}
0 & a_{12} & 0 \\
a_{21} & a_{22} & 0 \\
0 & a_{32} & a_{33} \\  
\end{array}
\right)a_{33}, a_{32}, a_{21}\neq 0, a_{12}= -a_{21},$ \cr
&$\mu_{3}=0, \mu_{2}=-1, \alpha= \dfrac{a_{21}\mu_{1}}{a_{12}} $ \cr
&$ \left(
\begin{array}{ccc}
0 & 0 & 0 \\
0 & a_{22} & 0 \\
a_{31} & 0 & a_{33} \\  
\end{array}
\right)a_{33}, a_{31}\neq 0, \alpha= -\mu_{1} $\cr
&$ \left(
\begin{array}{ccc}
0 & 0 & 0 \\
0 & a_{22} & 0 \\
a_{31} & a_{32} & a_{33} \\  
\end{array}
\right)a_{33}, a_{32}, a_{31}\neq 0, \alpha= -\mu_{1}, \mu_{3}=\dfrac{a_{31}(\alpha-\mu_{2})}{a_{33}}$ \cr
&$ \left(
\begin{array}{ccc}
0 & 0 & 0 \\
a_{21} & a_{22} & 0 \\
a_{31} & 0 & a_{33} \\  
\end{array}
\right)a_{33}, a_{31}, a_{21}\neq 0, \alpha= -\mu_{1} $
\\ \hline 
\end{tabular}
\end{center}

\begin{center}
\begin{tabular}{c c c}
\hline
 Dendriform  & Solutions of
 \cr
 algebra     &   $D$-equation \cr
 $ \mathcal{A} $    &
    \\ \hline
 $ D^{1}_{\mu_{1}, \mu_{2}, \mu_{3}} $:
&$ \left(
\begin{array}{ccc}
0 & 0 & 0 \\
a_{21} & a_{22} & 0 \\
a_{31} & a_{32} & a_{33} \\  
\end{array}
\right)a_{33}, a_{32}, a_{31}, a_{21}\neq 0, \alpha= -\mu_{1}, \mu_{3}=0$ 
\cr
 &$ \left(
\begin{array}{ccc}
0 & 0 & 0 \\
0 & a_{22} & a_{23} \\
a_{31} & a_{32} & a_{33} \\  
\end{array}
\right)a_{33}, a_{31}, a_{23}\neq 0, \alpha= -\mu_{1}, \mu_{3}=\mu_{2}=0 $  \cr
&$ \left(
\begin{array}{ccc}
0 & 0 & 0 \\
0 & a_{22} & a_{23} \\
a_{31} & a_{32} & a_{33} \\  
\end{array}
\right)a_{33}, a_{31}, a_{23}\neq 0, \alpha= -\mu_{1}\neq 0, a_{23}=-a_{32} $ \cr
&$ \left(
\begin{array}{ccc}
0 & 0 & 0 \\
a_{21} & a_{22} & a_{23} \\
a_{31} & a_{32} & a_{33} \\  
\end{array}
\right)a_{33}, a_{31}, a_{32}, a_{21}, a_{23}\neq 0, \alpha= -\mu_{1},$ \cr
&$\mu_{1}=\mu_{2}=\mu_{3}=0   $ \cr
&$ \left(
\begin{array}{ccc}
0 & 0 & 0 \\
a_{21} & a_{22} & a_{23} \\
a_{31} & a_{32} & a_{33} \\  
\end{array}
\right)a_{33}, a_{31}, a_{32}, a_{21}, a_{23}\neq 0; \alpha= -\mu_{1}\neq 0, $ \cr
& $ \mu_{3}=0; a_{31}a_{23}=a_{33}a_{21} $\cr
& $ \left(
\begin{array}{ccc}
0 & 0 & 0 \\
a_{21} & a_{22} & a_{23} \\
a_{31} & 0 & a_{33} \\  
\end{array}
\right) a_{21}, a_{23}\neq 0,$ \cr
&$a_{31}a_{23}\mu_{1}+ a_{33}a_{21}\mu_{2} + a_{33}a_{23}\mu_{3}=0,$ \cr
&$a_{31}a_{23}\mu_{2}+ a_{33}(a_{23}\mu_{3} + a_{21}\mu_{1})=0$ \cr
&$ \left(
\begin{array}{ccc}
0 & a_{12} & 0 \\
a_{21} & a_{22} & 0 \\
a_{31} & 0 & a_{33} \\  
\end{array}
\right) a_{31}, a_{33}, a_{12}\neq 0,$ \cr
& $a_{12}= -\mu_{2}a_{21}, \alpha=\mu_{1}=0$  \cr 
&$ \left(
\begin{array}{ccc}
0 & a_{12} & 0 \\
a_{21} & a_{22} & a_{23} \\
a_{31} & a_{32} & a_{33} \\  
\end{array}
\right) a_{31}, a_{33}, a_{12}, a_{23}\neq 0, \alpha=\mu_{1}=0,$\cr
&$a_{12}=-a_{21}\mu_{2}-a_{23}\mu_{3}, \mu_{2}=\dfrac{-a_{33}\mu_{3}}{a_{31}} $ \cr
& $ \left(
\begin{array}{ccc}
0 & a_{12} & a_{13} \\
0 & a_{22} & a_{23} \\
0 & 0 & 0 \\  
\end{array}
\right)a_{13}\neq 0, \mu_{1}=1, \mu_{3}= 0 $\cr
&$ \left(
\begin{array}{ccc}
0 & a_{12} & a_{13} \\
0 & a_{22} & a_{23} \\
0 & a_{32} & a_{33} \\  
\end{array}
\right)a_{13}, a_{33}\neq 0, a_{12}=\dfrac{a_{32}a_{13}}{a_{33}},$ \cr
&$ \alpha=-1, \mu_{1}=1, \mu_{3}=0$ \cr
&$ \left(
\begin{array}{ccc}
a_{11} & a_{12} & 0 \\
a_{21} & a_{22} & 0 \\
0 & 0 & 0 \\  
\end{array}
\right) a_{11}\neq 0$\cr
&$  \left(
\begin{array}{ccc}
a_{11} & a_{12} & 0 \\
a_{21} & a_{22} & 0 \\
0 & a_{32} & 0 \\  
\end{array}
\right)a_{11}, a_{32}\neq 0, \alpha= \mu_{2}=0$ \cr
&$ \left(
\begin{array}{ccc}
a_{11} & a_{12} & 0 \\
a_{21} & a_{22} & a_{23} \\
0 & a_{32} & 0 \\  
\end{array}
\right)a_{11}, a_{23}\neq 0, \mu_{1}=1-\mu_{2},$ \cr
& $ \mu_{2}=\dfrac{\alpha a_{32}}{a_{23}}, a_{12}=\dfrac{a^{2}_{23}\mu_{3}}{a_{32}}- \alpha a_{21} $  
\\ \hline
\end{tabular}
\end{center} 

\begin{center}
\begin{tabular}{c c c}
\hline
 Dendriform  & Solutions of
 \cr
 algebra     &   $D$-equation \cr
 $ \mathcal{A} $    &
    \\ \hline
 $ D^{1}_{\mu_{1}, \mu_{2}, \mu_{3}} $:& $\left(
\begin{array}{ccc}
a_{11} & a_{12} & a_{13} \\
0 & a_{22} & 0 \\
0 & 0 & 0 \\  
\end{array}
\right)a_{11}, a_{13}\neq 0,$ \cr
&$ \mu_{1}=1, \mu_{3}=\dfrac{-a_{11}(1+\mu_{2})}{a_{13}}$ \cr
&$ \left(
\begin{array}{ccc}
a_{11} & a_{12} & a_{13} \\
0 & a_{22} & a_{23} \\
0 & 0 & 0 \\  
\end{array}
\right)a_{11}, a_{13}, a_{23}\neq 0, \mu_{1}=1,$ \cr
&$ \mu_{2}=-1, \mu_{3}= 0, a_{21}=\dfrac{a_{23}a_{11}}{a_{13}}$ \cr
&$ \left(
\begin{array}{ccc}
a_{11} & a_{12} & a_{13} \\
-1 & a_{22} & a_{23} \\
0 & 0 & 0 \\  
\end{array}
\right)a_{11}, a_{13}, a_{23}\neq 0,$ \cr
&$\mu_{1}=1, \mu_{2}\neq -1, \mu_{3}= \dfrac{-a_{11}(1+\mu_{2})}{a_{13}}$ \cr
&$ \left(
\begin{array}{ccc}
a_{11} & a_{12} & 0 \\
a_{21} & a_{22} & 0 \\
a_{31} & a_{32} & 0 \\  
\end{array}
\right)a_{11}, a_{23}\neq 0,$\cr
&$ \alpha=\mu_{2}= -1, a_{12}= \dfrac{a_{32}a_{11}}{a_{31}}$ 
\\ \hline
$ D^{2}_{\mu_{1}, \mu_{2}}; \mu_{1}\neq 0$: & \cr
$ e_{3}\prec e_{3}=\mu_{1}e_{1}+\mu_{2}e_{2} $&  $ \left(
\begin{array}{ccc}
a_{11} & a_{12} & a_{13} \\
a_{21} & a_{22} & 0 \\
0 & 0 & 0 \\  
\end{array}
\right);  \left(
\begin{array}{ccc}
a_{11} & 0 & a_{13} \\
a_{21} & a_{22} & 0 \\
0 & a_{32} & 0 \\  
\end{array}
\right) a_{32}\neq 0$\cr
$e_{3}\prec e_{1}=\alpha e_{2}$ & $ \left(
\begin{array}{ccc}
a_{11} & a_{12} & a_{13} \\
a_{21} & a_{22} & 0 \\
0 & a_{32} & 0 \\  
\end{array}
\right)a_{32}, a_{12}\neq 0, \alpha= -1$\cr
$ e_{3}\succ e_{3}=-\mu_{1}e_{1}-\mu_{2}e_{2}$& $\left(
\begin{array}{ccc}
a_{11} & a_{12} & a_{13} \\
a_{21} & a_{22} & 0 \\
0 & a_{32} & 0 \\  
\end{array}
\right)a_{23}, a_{11}\neq 0, \alpha= 0,$\cr
$e_{1}\succ e_{3}= e_{2}$& $ a_{11}= -a_{23}\mu_{1}, a_{12}=- a_{23}\mu_{2}$ \cr
&$ \left(
\begin{array}{ccc}
a_{11} & a_{12} & a_{13} \\
a_{21} & a_{22} & a_{23} \\
0 & a_{32} & 0 \\  
\end{array}
\right)a_{23}, a_{11}, a_{32}\neq 0,$  \cr 
&$ \alpha=\dfrac{a^{2}_{23}\mu_{1} + a_{11}a_{23}}{a_{32}a_{11}},$ \cr
&$ a_{21}=\dfrac{a_{12}a_{32}(\alpha + 1)- a^{2}_{23}\mu_{2}- a_{12}a_{23}+ \mu_{2}a_{32}a_{23}}{\alpha a_{23}}  $ \cr
&$\left(
\begin{array}{ccc}
0 & 0 & a_{13} \\
a_{21} & a_{22} & 0 \\
a_{31} & a_{32} & 0 \\  
\end{array}
\right)a_{31}\neq 0$ \cr
&$\left(
\begin{array}{ccc}
a_{11} & a_{12} & a_{13} \\
a_{21} & a_{22} & 0 \\
a_{31} & a_{32} & 0 \\  
\end{array}
\right)a_{11}, a_{31}\neq 0, \alpha=-1,$ \cr
&$ a_{12}=\dfrac{a_{11}a_{32}}{a_{31}} $ \cr
& $\left(
\begin{array}{ccc}
a_{11} & a_{12} & a_{13} \\
a_{21} & a_{22} & a_{23} \\
a_{31} & a_{32} & 0 \\  
\end{array}
\right)a_{31}, a_{23}\neq 0, \alpha=0,$ \cr
&$a_{11}= -a_{23}\mu_{1}, a_{12}=-\mu_{2}a_{23}$ 
\\ \hline
\end{tabular}
\end{center}

\begin{center}
\begin{tabular}{c c c}
\hline
 Dendriform  & Solutions of
 \cr
 algebra     &   $D$-equation \cr
 $ \mathcal{A} $    & 
 \\ \hline
$ D^{3}_{\mu_{1}, \mu_{2}, \mu_{3}}; \mu_{1}\neq 0$: & \cr
$ e_{1}\prec e_{1}=\mu_{1}e_{2}  $ & $\left(
\begin{array}{ccc}
0 & a_{12} & 0 \\
a_{21} & a_{22} & a_{23} \\
0 & a_{32} & 0 \\  
\end{array}
\right); \left(
\begin{array}{ccc}
0 & 0 & 0 \\
0 & a_{22} & a_{23} \\
0 & a_{32} & a_{33} \\  
\end{array}
\right)a_{33}\neq 0, \mu_{3}= 0$     \cr
$ e_{3}\prec e_{1}=\mu_{2}e_{2}  $ & $ \left(
\begin{array}{ccc}
0 & 0 & 0 \\
a_{21} & a_{22} & a_{23} \\
0 & a_{32} & 0 \\  
\end{array}
\right)a_{33}, a_{21}\neq 0, \mu_{3}=0, \mu_{1}=\dfrac{-a_{23}\mu_{2}}{a_{21}}  $  \cr
$ e_{3}\prec e_{3}=\mu_{3}e_{2}  $ &  $ \left(
\begin{array}{ccc}
0 & 1 & 0\\
1 & a_{22} & 0 \\
0 & 0 & a_{33} \\  
\end{array}
\right)a_{33},a_{12}\neq 0, \alpha= 0, a_{12}=-\mu_{2}$   \cr
$ e_{1}\succ e_{1}=-\mu_{1}e_{2}  $ &  $\left(
\begin{array}{ccc}
0 & 1 & 0 \\
1 & a_{22} & a_{23} \\
0 & a_{32} & a_{33} \\  
\end{array}
\right)a_{33}, a_{23}\neq 0, \alpha=0, \mu_{2}=-1$    \cr
$ e_{3}\succ e_{1}=(\alpha-\mu_{2})e_{2}  $ & $ \left(
\begin{array}{ccc}
0 & a_{12} & 0 \\
1 & a_{22} & a_{23} \\
0 & a_{32} & a_{33} \\  
\end{array}
\right)a_{12}\neq 0,1, \alpha=0, \mu_{1}=\dfrac{a_{23}\mu_{2}+ a_{12}a_{23}}{a_{12}-1} $   \cr
$ e_{3}\succ e_{3}=-\mu_{3}e_{2}$ &  $\left(
\begin{array}{ccc}
0 & 0 & 0 \\
0 & a_{22} & 0 \\
a_{31} & a_{32} & a_{33} \\  
\end{array}
\right) a_{31}, a_{33}\neq 0, \alpha=\mu_{2}=\dfrac{-a_{31}\mu_{1}}{a_{33}}, \mu_{3}=0$    \cr
$ e_{1}\succ e_{3}= e_{2}  $ &
$\left(
\begin{array}{ccc}
0 & 0 & 0 \\
a_{21} & a_{22} & a_{23}\\
a_{31} & a_{32} & a_{33} \\  
\end{array}
\right)a_{31}, a_{33}, a_{21}\neq 0, \mu_{3}= 0,$ \cr
&$\alpha=\mu_{2}= \dfrac{-a_{31}\mu_{1}}{a_{33}}, a_{23}=\dfrac{-a_{21}\mu_{1}}{a_{33}}   $  \cr 
&    $ \left(
\begin{array}{ccc}
0 & 0 & 0 \\
1 & a_{22} & a_{23}\\
a_{31} & 0 & a_{33} \\  
\end{array}
\right)a_{31}, a_{33}, a_{23}, \mu_{3}\neq 0, \alpha\neq \mu_{2}, $   \cr 
& $a_{23}=\dfrac{(\alpha-\mu_{2})}{\mu_{3}}, \mu_{3}=\dfrac{a_{31}(\alpha-\mu_{2})}{a_{33}} $ \cr
&$\alpha=\dfrac{a_{31}\mu_{1}}{a_{33}}, \mu_{1}= \dfrac{-a_{33}(\mu_{2} + a_{23}\mu_{3})}{a_{31}}  $ \cr
 &$ \left(
\begin{array}{ccc}
a_{11} & a_{12} & a_{13} \\
0 & a_{22} & 0\\
0 & 0 & 0 \\  
\end{array}
\right)a_{11}, a_{13}\neq 0, a_{13}= a_{11}\mu_{1},$ \cr
& $ \mu_{3}=\dfrac{a_{11}(1+\mu_{2})}{a_{13}} $ \cr
& $\left(
\begin{array}{ccc}
a_{11} & a_{12} & a_{13} \\
a_{21} & a_{22} & a_{23}\\
0 & 0 & 0 \\  
\end{array}
\right)a_{11}, a_{13}, a_{21}\neq 0, a_{23}=a_{21}\mu_{1},$ \cr
& $ a_{13}=a_{11}\mu_{1}, 1 + \mu_{2} + \mu_{1}\mu_{3}= 0 $   \cr
& $\left(
\begin{array}{ccc}
a_{11} & a_{12} & a_{13} \\
a_{21} & a_{22} & a_{23}\\
a_{31} & a_{32} & a_{33} \\  
\end{array}
\right)a_{11}, a_{31}, a_{13}\neq 0, \alpha= -1, $  \cr
& $a_{33}=\dfrac{a_{31}a_{13}}{a_{11}}, a^{2}_{11}\mu_{1}+ a_{11}a_{13}\mu_{2}+ a^{2}_{13}\mu_{3}=0$  
\\ \hline 
\end{tabular}
\end{center}  

\begin{center}
\begin{tabular}{c c c}
\hline
 Dendriform  & Solutions of
 \cr
 algebra     &   $D$-equation 
\\ \hline
 $ D^{4}_{\mu_{1}, \mu_{2}}; \mu_{1}\neq 0$: & $ \left(
\begin{array}{ccc}
0 & 0 & 0 \\
0 & a_{22} & 0\\
0 & a_{32} & 0 \\  
\end{array}
\right); \left(
\begin{array}{ccc}
0 & 0 & 0 \\
0 & a_{22} & a_{23}\\
0 & a_{32} & 0 \\  
\end{array}
\right) a_{23}\neq 0, \mu_{2}=0$\cr
$ e_{1}\prec e_{3}= e_{2}  $ & $\left(
\begin{array}{ccc}
0 & 0 & 0 \\
0 & a_{22} & a_{23}\\
0 & a_{32} & 0 \\  
\end{array}
\right)a_{23}, \mu_{2}\neq 0, a_{23}= a_{32}$ \cr 
$ e_{3}\prec e_{3}= \mu_{2}e_{2}$ & $\left(
\begin{array}{ccc}
0 & 0 & 0 \\
0 & a_{22} & a_{23}\\
0 & a_{32} & a_{33} \\  
\end{array}
\right)a_{33}\neq 0, \mu_{2}=0$ \cr
$ e_{1}\prec e_{1}= \mu_{1}e_{3}  $ & $\left(
\begin{array}{ccc}
0 & a_{12} & 0 \\
0 & a_{22} & a_{23}\\
0 & 0 & 0 \\  
\end{array}
\right)a_{12}\neq 0, \mu_{2}= 0$\cr
$ e_{1}\prec e_{2}= \mu_{1}\mu_{2}e_{2}  $ &$\left(
\begin{array}{ccc}
0 & a_{12} & 0 \\
0 & a_{22} & a_{23}\\
0 & a_{32} & 0 \\  
\end{array}
\right)a_{12}, a_{32}\neq 0, \mu_{2}= \alpha=0 $\cr
$ e_{2}\prec e_{1}= \alpha\mu_{1}\mu_{2}e_{2}$& $\left(
\begin{array}{ccc}
0 & a_{12} & a_{13} \\
0 & a_{22} & a_{23}\\
0 & 0 & 0 \\  
\end{array}
\right)a_{13}\neq 0, \mu_{2}=0$\cr
$ e_{3}\succ e_{1}= \alpha e_{2}  $& $\left(
\begin{array}{ccc}
0 & a_{12} & 0 \\
a_{21} & a_{22} & a_{23}\\
0 & a_{32} & a_{33} \\  
\end{array}
\right)a_{21}, a_{33}\neq 0, \mu_{2}=0$\cr
$ e_{3}\succ e_{3}= -\mu_{2} e_{2}  $& $ a_{12}=\alpha a_{21} + \dfrac{a^{2}_{21}\mu_{1}}{a_{33}}, $\cr
$ e_{1}\succ e_{1}= -\mu_{1} e_{3}  $ & $ a_{23}= \dfrac{\alpha a_{32}(a_{12}-a_{21})}{a_{21}}+ $\cr
$ e_{1}\succ e_{2}= -\mu_{1}\mu_{2} e_{2}$& $ (\alpha^{2}-1)a^{2}_{33}+ 2\alpha a_{33}a_{21}\mu_{1}+ a^{2}_{21}\mu^{2}_{1} $\cr
$ e_{2}\succ e_{1}= -\alpha\mu_{1} \mu_{2} e_{2}$ &
 \\ \hline
\end{tabular}
\end{center}

{\bf Table 6}:
  Symmetric solutions and Heisenberg Connes cocycle structures.
\begin{center}
\begin{tabular}{c c c}
\hline
Dendriform   & Symmetric  &  Connes cocycles \cr
 algebra   & solutions     &  structures over $ \mathcal{H}\oplus  \mathcal{H}^{\ast} $\cr
structures &&
\\ \hline
 $ D^{1}_{\mu_{1}, \mu_{2}, \mu_{3}}$ &$ a_{22}e_{2}\otimes e_{2}  $ &
                        $e_{3}\ast e^{\ast}_{2}= -a_{11}\mu_{1} e_{1};$ $  e^{\ast}_{3}\ast e^{\ast}_{2}= -a_{33}\mu_{3} e^{\ast}_{3} $  \cr
\cline{2-2}
    & $ a_{22}e_{2}\otimes e_{2} +  $  & $ e^{\ast}_{1}\ast e^{\ast}_{2}= -a_{11}\mu_{2} e^{\ast}_{3}$; $ e^{\ast}_{1}\ast e^{\ast}_{1}=(\mu_{1}-1)a_{11} e^{\ast}_{3}$\cr
   &   $ a_{23}(e_{2}\otimes e_{3} + e_{3}\otimes e_{2}) $ &
                         $e^{\ast}_{2}\ast e^{\ast}_{2}=(a_{23}\mu_{2}-a_{23}\mu_{1}-\alpha a_{23}) e^{\ast}_{1} +$ \cr
   &$a_{23}\neq 0$ &                      $ (-a_{12}\mu_{2}+a_{12}\mu_{1}-\alpha a_{23}- a_{12}) e^{\ast}_{3}$ \cr 
 & $ \mu_{3}=\dfrac{-a_{21}(\mu_{1}+ \mu_{2})+ a_{23}\mu_{1}}{a_{23}}$  & $ e^{\ast}_{2}\ast e^{\ast}_{3}=(a_{33}\mu_{2}-a_{33}\alpha)e^{\ast}_{1} + a_{33}\mu_{3} e^{\ast}_{3} $  \cr
  \cline{2-2}
 & $ a_{12}(e_{1}\otimes e_{2} + e_{2}\otimes e_{1})$   &  $ e_{1}\ast e^{\ast}_{2}= (-1+\mu_{2})a_{23}e_{2} + a_{33}\mu_{2}e_{3} + \mu_{2}e^{\ast}_{3}$ \cr
 &    $a_{22}e_{2}\otimes e_{2}; a_{12}\neq 0$    & $ e^{\ast}_{2}\ast e_{1}=  (1-\mu_{1})e^{\ast}_{1}+ (1-\mu_{1})a_{33}e_{3} +$   \cr
 \cline{2-2}
& $ a_{12}(e_{1}\otimes e_{2} + e_{2}\otimes e_{1})+$ &$(1-\mu_{1}-\alpha)a_{23}e_{2} $ \cr
&$ a_{23}(e_{2}\otimes e_{3} + e_{3}\otimes e_{2})+$ &$ e_{3}\ast e^{\ast}_{2}=(-a_{12}\alpha +a_{12}\mu_{1} + a_{23}\mu_{3})e_{2}+ $ \cr
& $ a_{22}e_{2}\otimes e_{2}; a_{12}, a_{23}\neq 0$& $\mu_{1}e^{\ast}_{1} + \mu_{3}e^{\ast}_{3} +a_{33}\mu_{3}e_{3}$\cr
&$ \mu_{3}=\dfrac{a_{12}}{a_{23}}-a_{22}+ \mu_{1} $& $e^{\ast}_{2}\ast e_{3}=a_{11}(\alpha- \mu_{2})e_{1}- \mu_{3}a_{33}e_{3}+$ \cr
\cline{2-2}
&$a_{23}(e_{2}\otimes e_{3} + e_{3}\otimes e_{2})+$& $(a_{12}(\alpha- \mu_{2}) -a_{23}\mu_{3}- a_{12} -\alpha a_{23})e_{2}  $\cr
&$a_{22}e_{2}\otimes e_{2}; a_{21}=a_{31}=0$& $ e^{\ast}_{3}\ast e^{\ast}_{2}= -a_{33}\mu_{1}e^{\ast}_{1} $\cr
&$\mu_{1}=\mu_{2}=\mu_{3}=0$&$ e_{1}\ast e_{3}= e_{2}$; $ e_{3}\ast e_{1}=\alpha e_{2}$\cr
\cline{2-2}
& $a_{23}(e_{2}\otimes e_{3} + e_{3}\otimes e_{2})$ &$e^{\ast}_{3}\ast e_{1}= -a_{33}\alpha e_{2};$ $e_{3}\ast e^{\ast}_{1}= -a_{11}\alpha e_{2}$\cr
& $ a_{22}e_{2}\otimes e_{2}; a_{33}\neq 0; \mu_{3}=0 $ & $e_{1}\ast e^{\ast}_{3}= (-a_{33} -a_{11})e_{2}$  
\\ \hline
\end{tabular}
\end{center}

\begin{center}
\begin{tabular}{c c c}
\hline
Dendriform   & Symmetric  &  Connes cocycles \cr
 algebra   & solutions     &  structures over $ \mathcal{H}\oplus  \mathcal{H}^{\ast} $\cr
structures &&
\\ \hline
$D^{1}_{\mu_{1}, \mu_{2}, \mu_{3}}$ & $ a_{23}(e_{2}\otimes e_{3} + e_{3}\otimes e_{2})+$& $e_{3}\ast e^{\ast}_{2}= -a_{11}\mu_{1} e_{1}$  \cr
&$ a_{12}(e_{1}\otimes e_{2} + e_{2}\otimes e_{1})$+&$ e^{\ast}_{1}\ast e^{\ast}_{2}= -a_{11}\mu_{2} e^{\ast}_{3}$  \cr
&$a_{22}e_{2}\otimes e_{2} + a_{33}e_{3}\otimes e_{3}  $& $e_{1}\ast e^{\ast}_{3}= (-a_{33} -a_{11})e_{2}$  \cr
&$\mu_{3}=0; \alpha=\mu_{1}=1-\mu_{2}$& $ e^{\ast}_{1}\ast e^{\ast}_{1}=(\mu_{1}-1)a_{11} e^{\ast}_{3}$ \cr
\cline{2-2}
&$ a_{12}(e_{1}\otimes e_{2} + e_{2}\otimes e_{1})$+& $ e^{\ast}_{2}\ast e^{\ast}_{3}=(a_{33}\mu_{2}-a_{33}\alpha)e^{\ast}_{1}$ \cr
&$a_{22}e_{2}\otimes e_{2} + a_{33}e_{3}\otimes e_{3}$& $+ a_{33}\mu_{3} e^{\ast}_{3}$ \cr
&$a_{33}, a_{21}, a_{12}\neq 0$& $ e^{\ast}_{2}\ast e_{1}=  (1-\mu_{1})e^{\ast}_{1}+ $  \cr
&$\mu_{3}=0; \alpha= \mu_{1}$ & $(1-\mu_{1}-\alpha)a_{23}e_{2}+ $\cr
\cline{2-2}
&$ a_{12}(e_{1}\otimes e_{2} + e_{2}\otimes e_{1})$+& $(1-\mu_{1})a_{33}e_{3}$  \cr
&$ a_{22}e_{2}\otimes e_{2} + a_{11}e_{1}\otimes e_{1} $&$ e_{3}\ast e^{\ast}_{2}=(-a_{12}\alpha + $\cr
&$a_{11}\neq 0$&$a_{12}\mu_{1} + a_{23}\mu_{3})e_{2}+$\cr
\cline{2-2}
&$ a_{12}(e_{1}\otimes e_{2} + e_{2}\otimes e_{1})$+&$\mu_{1}e^{\ast}_{1} + \mu_{3}e^{\ast}_{3} +a_{33}\mu_{3}e_{3}$ \cr
&$ a_{22}e_{2}\otimes e_{2} + a_{11}e_{1}\otimes e_{1} $&$ e^{\ast}_{3}\ast e^{\ast}_{2}= -a_{33}\mu_{1}e^{\ast}_{1}$\cr
&$a_{23}(e_{2}\otimes e_{3} + e_{3}\otimes e_{2})$&$ e_{1}\ast e_{3}= e_{2}$; $ e_{3}\ast e_{1}=\alpha e_{2}$\cr
&$\mu_{1}= 1-\alpha=\mu_{2}$&$e^{\ast}_{3}\ast e_{1}= -a_{33}\alpha e_{2}$\cr
&$a_{12}=\dfrac{a_{23}\mu_{3}}{1+\alpha}; \alpha\neq -1$&$e_{3}\ast e^{\ast}_{1}= -a_{11}\alpha e_{2}$\cr
&&$ e_{1}\ast e^{\ast}_{2}= (-1+\mu_{2})a_{23}e_{2}$\cr
&&$+ a_{33}\mu_{2}e_{3} + \mu_{2}e^{\ast}_{3}$ \cr
&&$  e^{\ast}_{3}\ast e^{\ast}_{2}= -a_{33}\mu_{3} e^{\ast}_{3} $\cr
&&$e^{\ast}_{2}\ast e^{\ast}_{2}=(a_{23}\mu_{2}$ \cr
&&$-a_{23}\mu_{1}-\alpha a_{23}) e^{\ast}_{1} +$\cr
&& $ (-a_{12}\mu_{2}+a_{12}\mu_{1}$\cr
&&$-\alpha a_{23}- a_{12})e^{\ast}_{3}$\cr
&&$e^{\ast}_{2}\ast e_{3}=a_{11}(\alpha- \mu_{2})e_{1}- \mu_{3}a_{33}e_{3}+$\cr
&&$(a_{12}(\alpha- \mu_{2}) -a_{23}\mu_{3}- a_{12} -\alpha a_{23})e_{2}  $
\\ \hline
$ D^{2}_{\mu_{1}, \mu_{2}}; \mu_{1}\neq 0$: &$a_{12}(e_{1}\otimes e_{2} + e_{2}\otimes e_{1})  $ &
$e^{\ast}_{2}\ast e^{\ast}_{1}=(a_{13}\mu_{2}-a_{13}-\mu_{1}a_{23}$ \cr
&$ a_{11}e_{1}\otimes e_{1} + a_{22}e_{2}\otimes e_{2}$  & $-a_{11}-\alpha a_{13})e^{\ast}_{3}+$   \cr
\cline{2-2}
    & $ a_{22}e_{2}\otimes e_{2} + a_{11}e_{1}\otimes e_{1}$& $  (-a_{23}+a_{13}\alpha)e^{\ast}_{1} $\cr
   &   $ a_{12}(e_{1}\otimes e_{2} + e_{2}\otimes e_{1}) $ &
                         $  e^{\ast}_{2}\ast e^{\ast}_{2}=\alpha a_{23} e^{\ast}_{1} +(\alpha- 1)a_{12}e^{\ast}_{3}$ \cr
   &$ a_{23}(e_{2}\otimes e_{3} + e_{3}\otimes e_{2}) $ &$e_{1}\ast e^{\ast}_{2}=(-a_{23}+ \alpha a_{13})e_{1}+ $ \cr 
 & $\alpha=\dfrac{a_{23}\mu_{1}+ a_{11}}{a_{11}}$  &  $\alpha a_{23}e_{2} + \alpha e^{\ast}_{3}$ \cr
  
 & $a_{23}, a_{11}, a_{32}\neq 0$   &  $ e_{3}\ast e^{\ast}_{3}= -a_{13}\alpha e_{2}$ \cr
 \cline{2-2}
 &    $ a_{13}(e_{1}\otimes e_{3} + e_{3}\otimes e_{1})+$&  $e_{1}\ast e^{\ast}_{1}= -a_{13}e_{2}$   \cr
& $ a_{22}e_{2}\otimes e_{2}; a_{13}\neq 0 $ & $ e_{3}\ast e^{\ast}_{1}=(-a_{11}\alpha +a_{23}\mu_{1})e_{2}$ \cr
  \cline{2-2}
&$ a_{23}(e_{2}\otimes e_{3} + e_{3}\otimes e_{2})+$ &$ + a_{13}\mu_{1}e_{1} + \mu_{1}e^{\ast}_{3} $ \cr
& $ a_{22}e_{2}\otimes e_{2}+ a_{11}e_{1}\otimes e_{1}$& $e^{\ast}_{1}\ast e_{3}= (-\mu_{1}a_{23}-a_{11}-\alpha a_{13})e_{2}$\cr
&$a_{12}(e_{1}\otimes e_{2} + e_{2}\otimes e_{1})$& $ e^{\ast}_{3}\ast e_{3}=-a_{13}e_{2}$ \cr
&$a_{13}(e_{1}\otimes e_{3} + e_{3}\otimes e_{1})+$&$ e_{1}\ast e_{3}= e_{2}$\cr
&$a_{31}, a_{23}\neq 0; \alpha= 0$& $ e_{3}\ast e_{1}=\alpha e_{2}$ \cr
&$a_{11}=-a_{23}\mu_{1}$& $ e^{\ast}_{2}\ast e^{\ast}_{3}=-a_{13}e^{\ast}_{3}$ \cr
&$ a_{12}=-\mu_{2}a_{23}$& $ e^{\ast}_{3}\ast e^{\ast}_{2}= -a_{13}\alpha e^{\ast}_{3}$\cr
&&$e_{3}\ast e^{\ast}_{2}=(a_{12}\alpha+a_{23}\mu_{2})e_{2} +$\cr
&& $ a_{13}\mu_{2} e_{1} + \mu_{2}e^{\ast}_{3} $\cr
&& $ e^{\ast}_{2}\ast e_{1}= a_{13}e_{1}+ e^{\ast}_{3}$\cr
&& $ e^{\ast}_{2}\ast e_{3}=-\mu_{2}a_{13}e_{1}+ $\cr
&&$(-\mu_{2}a_{23}-a_{12})e_{2}- \mu_{2}e^{\ast}_{3}$\cr
&&$ e^{\ast}_{1}\ast e_{3}=-a_{13}e_{2}$\cr
&&$ -\mu_{1}a_{13}e_{1} -\mu_{1}e_{3}$\cr
&&$e^{\ast}_{1}\ast e^{\ast}_{2}= (-a_{11}\alpha$\cr
&&$ + a_{23}\mu_{1}-\mu_{2}a_{13})e^{\ast}_{3}$ 
\\ \hline
\end{tabular}
\end{center}

\begin{center}
\begin{tabular}{c c c}
\hline
Dendriform   & Symmetric  &  Connes cocycles \cr
 algebra   & solutions     &  structures over $ \mathcal{H}\oplus  \mathcal{H}^{\ast} $\cr
structures &&
\\ \hline
$ D^{3}_{\mu_{1}, \mu_{2}, \mu_{3}}; \mu_{1}\neq 0$: & $a_{12}(e_{1}\otimes e_{2} + e_{2}\otimes e_{1})  $  & $ e_{1}\ast e_{3}= e_{2}$; $ e_{3}\ast e_{1}=\alpha e_{2};$ $e_{1}\ast e^{\ast}_{1}= -a_{13}e_{2}$ \cr
                      & $ a_{23}(e_{2}\otimes e_{3} + e_{3}\otimes e_{2})$  &  $e^{\ast}_{1}\ast e^{\ast}_{2}= -a_{11}\mu_{1}e^{\ast}_{1} + (-a_{11}\mu_{2}-a_{13}\mu_{3})e^{\ast}_{3}$\cr
                        &$a_{22}(e_{2}\otimes e_{2}) $&
$e^{\ast}_{2}\ast e^{\ast}_{1}=(a_{11}\mu_{1}+a_{13}\mu_{2}-\alpha a_{13})e^{\ast}_{1}+$ \cr
&&  $  (a_{13}\mu_{3}- a_{11})e^{\ast}_{3} $  \cr
\cline{2-2}
    & $ a_{22}e_{2}\otimes e_{2} + a_{33}e_{3}\otimes e_{3} $  & $ e^{\ast}_{3}\ast e^{\ast}_{2}= -\mu_{1}a_{13}e^{\ast}_{1} + (-a_{13}\mu_{2}-a_{33}\mu_{3})e^{\ast}_{3}$\cr
   &   $ a_{23}(e_{2}\otimes e_{3} + e_{3}\otimes e_{2}) $ &
                         $  e^{\ast}_{2}\ast e^{\ast}_{2}=-\alpha a_{23} e^{\ast}_{1} +(-1- \mu_{2})a_{12}e^{\ast}_{3}$ \cr
   &$ a_{33}\neq 0; \mu_{3}=0 $ &$e_{1}\ast e^{\ast}_{3}=-a_{33}e_{2}; e^{\ast}_{1}\ast e_{3}= -a_{11}e^{\ast}_{2}$ \cr
   \cline{2-2} 
 & $a_{22}e_{2}\otimes e_{2} + a_{33}e_{3}\otimes e_{3}$ & $ e^{\ast}_{1}\ast e_{1}= -a_{13}\alpha e_{2} $  \cr
  
 & $a_{12}(e_{1}\otimes e_{2} + e_{2}\otimes e_{1})$ &  $ e^{\ast}_{3}\ast e_{1}= -a_{33}\alpha e_{2}; e_{3}\ast e^{\ast}_{3}= -a_{13}\alpha e_{2}$ \cr
 
 &    $ \alpha=0$    & $ e_{3}\ast e^{\ast}_{1}=(-a_{11}  - a_{33})\alpha e_{2}$   \cr
& $a_{12}=-\mu_{2}\neq 0$ & $ e^{\ast}_{2}\ast
 e^{\ast}_{3}=(\mu_{1}a_{13}+\mu_{2}a_{33}-a_{33}\alpha)e^{\ast}_{1}+$\cr
  \cline{2-2}
&$ a_{23}(e_{2}\otimes e_{3} + e_{3}\otimes e_{2})+$ &$ e_{3}\ast e^{\ast}_{2}=a_{13}\mu_{3}e_{1}+ a_{33}\mu_{3})e_{3} + $\cr
& $a_{22}e_{2}\otimes e_{2}+ a_{11}e_{1}\otimes e_{1}$&$(-a_{12}\alpha+ a_{23}\mu_{3})e_{2} + \mu_{3}e^{\ast}_{3} $ \cr
& $ e_{1}\otimes e_{2} + e_{2}\otimes e_{1}$ & $e_{1}\ast e^{\ast}_{2}= (\mu_{1}a_{11}+ a_{13}\mu_{2})e_{1}+ \mu_{1}e^{\ast}_{1}$ \cr
&$\mu_{2}=-1; \alpha= 0$& $ (-a_{23}+ a_{12}\mu_{1})e_{2} + \mu_{2}e^{\ast}_{3}$\cr
\cline{2-2}
&$a_{22}e_{2}\otimes e_{2}+ a_{11}e_{1}\otimes e_{1}$& $(a_{13}\mu_{1}+ a_{33}\mu_{2})e_{3}  $\cr
&$a_{12}(e_{1}\otimes e_{2} + e_{2}\otimes e_{1})$& $ e^{\ast}_{2}\ast e_{1}=(-a_{11}\mu_{1}+ a_{13})e_{1} -\mu_{1}e^{\ast}_{1}$ \cr
&$a_{13}(e_{1}\otimes e_{3} + e_{3}\otimes e_{1})$&$  (-a_{12}\mu_{1}+ a_{23} -a_{23}\alpha)e_{2}$ \cr
&$a_{33}e_{3}\otimes e_{3}$&$ (-\mu_{1}a_{13}+ a_{33})e_{3} +e^{\ast}_{3} $ \cr
&$a_{23}(e_{2}\otimes e_{3} + e_{3}\otimes e_{2})+$& $e^{\ast}_{2}\ast e_{3}=(\alpha-\mu_{2})e^{\ast}_{1} -\mu_{3}e^{\ast}_{3}+$ \cr
&$a_{31}=\dfrac{a^{2}_{13}}{a_{11}}$& $(a_{11}(\alpha- \mu_{2})-a_{13}\mu_{3})e_{1}+$\cr
&$a^{2}_{11}+ a_{11}a_{13}\mu_{2}+ a^{2}_{13}\mu_{3}=0$ 
\\ \hline
 $ D^{4}_{\mu_{1}, \mu_{2}}; \mu_{1}\neq 0$:   & $a_{22}e_{2}\otimes e_{2}$  & $ e_{1}\ast e_{3}= e_{2}$; $e_{3}\ast e_{1}=\alpha e_{2}$  \cr
 \cline{2-2}
& $ a_{23}(e_{2}\otimes e_{3} + e_{3}\otimes e_{2})$  &  $e_{2}\ast e^{\ast}_{2}= \mu_{1}\mu_{2}e^{\ast}_{1}$\cr
& $a_{22}e_{2}\otimes e_{2}$ &$e^{\ast}_{2}\ast e_{2}= -\alpha\mu_{1}\mu_{2}e^{\ast}_{1}$ \cr
 \cline{2-2}
&$a_{23}(e_{2}\otimes e_{3} + e_{3}\otimes e_{2})$& $e^{\ast}_{2}\ast e^{\ast}_{1}= a_{12}\alpha \mu_{1}\mu_{2}e^{\ast}_{1}$  \cr
    & $ a_{22}e_{2}\otimes e_{2} + a_{33}e_{3}\otimes e_{3} $  & $e^{\ast}_{1}\ast e^{\ast}_{2}= -a_{12}\mu_{1}\mu_{2}e^{\ast}_{1}$\cr
   &   $a_{33}\neq 0; \mu_{2}=0 $ &
                         $  e^{\ast}_{3}\ast e_{1}=(-\mu_{1}a_{12}-a_{33}\alpha)e_{2} -\mu_{1}e^{\ast}_{1}$ \cr
     \cline{2-2}
   &$ a_{23}(e_{2}\otimes e_{3} + e_{3}\otimes e_{2}) $ &$e_{1}\ast e^{\ast}_{3}=(-a_{33}+  a_{12}\mu_{1})e_{2} +\mu_{1}e^{\ast}_{1} $ \cr 
 & $a_{12}(e_{1}\otimes e_{2} + e_{2}\otimes e_{1}$  & $ e^{\ast}_{2}\ast e_{1}=-a_{12}\mu_{1}\mu_{2}e_{1} -a_{23}\mu_{1}\mu_{2} e_{3} -\mu_{1}\mu_{2}e^{\ast}_{2} $  \cr
  
 & $a_{22}e_{2}\otimes e_{2} + a_{33}e_{3}\otimes e_{3}$   &  $ e_{3}\ast e^{\ast}_{2}= (-a_{12}\alpha +a_{23}\mu_{2}+ a_{12})e_{2} + e^{\ast}_{1} +$ \cr
&    $ a_{21}, a_{33}\neq 0; \mu_{2}=0$    & $ a_{33}\mu_{2}e_{3} + \mu_{2}e^{\ast}_{3}+ \alpha e^{\ast}_{1}-\mu_{2}e^{\ast}_{3}$   \cr
& $ a_{12}= a_{21}(\alpha + \dfrac{a_{21}\mu_{1}}{a_{33}}) $ & $ e_{1}\ast
 e^{\ast}_{2}=(-1+\alpha\mu_{1}\mu_{2}+\mu_{2})a_{23}e_{2}+$ \cr
 &$a_{23}=\dfrac{\alpha a_{32}(a_{12}-a_{21})}{a_{21}}$& $(a_{23}\alpha\mu_{1}\mu_{2}+ a_{33}\mu_{2})e_{3}+$\cr
&$(\alpha^{2}-1)a^{2}_{33}+ a^{2}_{21}\mu^{2}_{1}+$ &$ a_{12}\alpha\mu_{1}\mu_{2}e_{1}+ \alpha \mu_{1}\mu_{2}e^{\ast}_{2}$ \cr
& $ 2a_{33}\alpha a_{21}\mu_{1}=0$& $e^{\ast}_{2}\ast e^{\ast}_{2}=(-a_{23}+ a_{22}\alpha\mu_{1}\mu_{2} +$\cr
&& $a_{23}\mu_{2}-a_{23}\alpha -a_{22}\mu_{1}\mu_{2})e^{\ast}_{1}$ \cr
&& $ e^{\ast}_{3}\ast e^{\ast}_{2}=(-a_{33} + a_{12}\mu_{1} -a_{23}\mu_{1}\mu_{2})e^{\ast}_{1}+$\cr
&& $-a_{33}\mu_{2}e^{\ast}_{3}$\cr
&&$e^{\ast}_{2}\ast e^{\ast}_{3}=(a_{23}\alpha \mu_{1}\mu_{2} + a_{33}\mu_{2}- \mu_{1}a_{12}$ \cr
&& $-a_{23}\alpha)e^{\ast}_{1}+ a_{33}\mu_{2}e^{\ast}_{3}$
 \\ \hline
\end{tabular}
\end{center}

\begin{center}
{\bf Table 7}:
 Solutions of the $ 3- $dimensional non decomposable associative Yang-Baxter equations
\begin{tabular}{c c}
\hline
Associative  & Solutions  \cr
 algebra $ \mathcal{A} $  & of the AYBE
 \\ \hline
$\mathcal{A}_{1}= \mathcal{H}_{1} $   & See Table 2.
\\ \hline
$\mathcal{A}_{2}= \mathcal{H}_{\alpha}; \alpha\neq 1 $   & See Table 2.
\\ \hline
$ \mathcal{A}_{3} $:$ e_{1}\cdot e_{1}= e_{2} $ & $ \left(
\begin{array}{ccc}
0 & 0 & 0 \\
0 & a_{22} & a_{23} \\
0 & a_{32} & a_{33} \\  
\end{array}
\right)
$;
 $ \left(
\begin{array}{ccc}
0 & 0 & a_{13} \\
0 & a_{22} & a_{23} \\
a_{31} & a_{32} & a_{33} \\  
\end{array}
\right) a_{22}, a_{13}, a_{31}\neq 0
,$\cr
$ e_{2}\cdot e_{1}= e_{3} $ & $ a_{32}=\frac{a_{23}(2a_{13}- a_{31})}{a_{13}- 2a_{31}} (a_{13}\neq 2a_{31})$   \cr       
$  e_{1}\cdot e_{2}= e_{3}  $ & $ \left(
\begin{array}{ccc}
0 & 0 & a_{13} \\
0 & a_{22} & 0 \\
a_{31} & a_{32} & a_{33} \\  
\end{array}
\right) a_{22}, a_{13}, a_{31}\neq 0, a_{13}= 2a_{31}
$ \cr
&  $ \left(
\begin{array}{ccc}
0 & 0 & a_{13} \\
0 & a_{22} & a_{23} \\
a_{31} & a_{32} & a_{33} \\  
\end{array}
\right) a_{22}, a_{23}, a_{13}, a_{31}\neq 0, a_{13}= 2a_{31}
 $ 
 \\ \hline
$ \mathcal{A}_{4} $: $ e_{1}\cdot e_{3}= e_{2} $ & $ \left(
\begin{array}{ccc}
a_{11} & a_{12} & 0 \\
a_{21} & a_{22} & 0 \\
0 & 0 & 0 \\  
\end{array}
\right)
  $;
  $ \left(
\begin{array}{ccc}
a_{11} & a_{12} & 0 \\
a_{21} & a_{22} & a_{23} \\
0 & 0 & 0 \\  
\end{array}
\right)  a_{21}= a_{12}, a_{23}\neq 0
 $ \cr
$ e_{3}\cdot e_{3}= e_{3}$ & $ \left(
\begin{array}{ccc}
a_{11} & a_{12} & 0 \\
a_{21} & a_{22} & a_{23} \\
0 & a_{32} & 0 \\  
\end{array}
\right) a_{23}= -a_{32}\neq 0, a_{21}= -a_{12}\neq 0$\cr
$ e_{2}\cdot e_{3}= e_{2} $ &
\\ \hline
$ \mathcal{A}_{5} $: $e_{2}\cdot e_{3}= e_{2}$ & $ \left(
\begin{array}{ccc}
a_{11} & a_{12} & 0 \\
a_{21} & a_{22} & 0 \\
0 & 0 & 0 \\  
\end{array}
\right) $;$ \left(
\begin{array}{ccc}
a_{11} & 0 & 0 \\
a_{21} & 0 & 0 \\
a_{31} & 0 & 0 \\  
\end{array}
\right)  a_{31}\neq 0$ \cr
$ e_{3}\cdot e_{1}= e_{1} $ & $ \left(
\begin{array}{ccc}
0 & a_{12} & a_{13} \\
a_{21} & 0 & 0 \\
a_{31} & 0 & 0 \\  
\end{array}
\right) a_{31}= -a_{13}\neq 0$
 \cr
 $ e_{3}\cdot e_{3}= e_{3} $  & $ \left(
\begin{array}{ccc}
0 & a_{12} & 0 \\
a_{21} & 0 & a_{23} \\
0 & a_{32} & 0 \\  
\end{array}
\right)a_{23}= -a_{32}\neq 0$\cr
&$ \left(
\begin{array}{ccc}
0 & 0 & 0 \\
a_{21} & a_{22} & a_{23} \\
0 & a_{32} & 0 \\  
\end{array}
\right) a_{23}\neq 0$ 
\\ \hline
$ \mathcal{A}_{6} $: $ e_{3}\cdot e_{1}= e_{2} $ & $ \left(
\begin{array}{ccc}
a_{11} & a_{12} & 0 \\
a_{21} & a_{22} & 0 \\
0 & 0 & 0 \\  
\end{array}
\right) a_{11}= -a_{12};$
 $ \left(
\begin{array}{ccc}
a_{11} & a_{12} & 0 \\
a_{21} & a_{22} & 0 \\
0 & 0 & 0 \\  
\end{array}
\right) a_{11}\neq -a_{12}$
 \cr
$ e_{3}\cdot e_{2}= e_{2}$ &$ \left(
\begin{array}{ccc}
a_{11} & a_{12} & 0 \\
0 & a_{22} & a_{23} \\
0 & a_{32} & 0 \\  
\end{array}
\right) a_{32}= -a_{23}\neq 0, a_{11}= -a_{12}; $ $\left(
\begin{array}{ccc}
0 & 0 & 0 \\
0 & a_{22} & 0 \\
0 & a_{32} & 0 \\  
\end{array}
\right)a_{32}\neq 0$\cr
$ e_{3}\cdot e_{3}= e_{3} $   &
\\  \hline
$ \mathcal{A}_{7} $: $ e_{1}\cdot e_{2}= e_{1} $ & 
 $ \left(
\begin{array}{ccc}
0 & 0 & 0 \\
0 & 0 & 0 \\
a_{31} & 0 & 0 \\  
\end{array}
\right)a_{31}\neq 0
$; $ \left(
\begin{array}{ccc}
0 & a_{12} & a_{13} \\
0 & 0 & 0 \\
0 & 0 & 0 \\  
\end{array}
\right)a_{13} =a_{12}\neq 0
$ \cr
$ e_{2}\cdot e_{2}= e_{2} $ &$ \left(
\begin{array}{ccc}
0 & 0 & a_{13} \\
0 & 0 & 0 \\
a_{31} & 0 & 0 \\  
\end{array}
\right)
 a_{31}= -a_{13}\neq 0$ \cr
 $ e_{3}\cdot e_{1}= e_{1} $ & $ \left(
\begin{array}{ccc}
a_{11} & a_{12} & 0 \\
0 & 0 & 0 \\
0 & 0 & 0 \\  
\end{array}
\right)
$; $ \left(
\begin{array}{ccc}
0 & 0 & 0 \\
a_{21} & 0 & 0 \\
a_{31} & 0 & 0 \\  
\end{array}
\right)  a_{31}= a_{21}\neq 0
$ \cr
$ e_{3}\cdot e_{3}= e_{3} $ & $ \left(
\begin{array}{ccc}
0 & a_{12} & 0 \\
a_{21} & 0 & 0 \\
0 & 0 & 0 \\  
\end{array}
\right)a_{21}= -a_{12} \neq 0
 $
\\ \hline 
\end{tabular}
\end{center}

\begin{center}
\begin{tabular}{c c}
\\ \hline
$ \mathcal{A}_{8} $: $ e_{1}\cdot e_{3}= e_{1}$ & $ \left(
\begin{array}{ccc}
a_{11} & a_{12} & 0 \\
a_{21} & a_{22} & 0 \\
0 &  & 0 \\  
\end{array}
\right) 
 $;
 $ \left(
\begin{array}{ccc}
a_{11} & a_{12} & a_{13} \\
0 & 0 & 0 \\
0 & 0 & 0 \\  
\end{array}
\right) a_{13}\neq 0;$ $ \left(
\begin{array}{ccc}
0 & 0 & 0 \\
a_{21} & 0 & 0 \\
a_{31} & 0 & 0 \\  
\end{array}
\right) a_{31} \neq 0 $\cr
$ e_{2}\cdot e_{3}= e_{2} $ & $ \left(
\begin{array}{ccc}
0 & 0 & 0 \\
a_{21} & a_{22} & a_{23} \\
0 & 0 & 0 \\  
\end{array}
\right)  a_{23}\neq 0;$ 
 $ \left(
\begin{array}{ccc}
0 & 0 & 0 \\
0 & 0 & a_{23} \\
0 & a_{32} & 0 \\  
\end{array}
\right)a_{32} =a_{23}\neq 0$ \cr
$ e_{3}\cdot e_{3}= e_{3} $ & \cr
$ e_{3}\cdot e_{1}= e_{1}  $ &
\\ \hline 
$ \mathcal{A}_{9} $: $ e_{2}\cdot e_{3}= e_{2}$ & $ \left(
\begin{array}{ccc}
a_{11} & a_{12} & 0 \\
a_{21} & a_{22} & 0 \\
0 & 0 & 0 \\  
\end{array}
\right) $;
 $ \left(
\begin{array}{ccc}
a_{11} & 0 & a_{13} \\
0 & 0 & a_{23} \\
a_{31} & 0 & 0 \\  
\end{array}
\right) a_{31}\neq 0 $; $ \left(
\begin{array}{ccc}
0 & 0 & 0 \\
a_{21} & 0 & a_{23} \\
0 & 0 & 0 \\  
\end{array}
\right)a_{23}\neq 0$ \cr
$ e_{3}\cdot e_{1}= e_{1} $ & \cr
$ e_{3}\cdot e_{2}= e_{2} $ &$ \left(
\begin{array}{ccc}
0 & 0 & a_{13} \\
0 & a_{22} & 0 \\
a_{31} & 0 & 0 \\  
\end{array}
\right) a_{31}= -a_{13}\neq 0 
$  \cr
$ e_{3}\cdot e_{3}= e_{3} $ &$ \left(
\begin{array}{ccc}
a_{11} & a_{12} & 0 \\
0 & 0 & 0 \\
a_{31} & a_{32} & 0 \\  
\end{array}
\right)a_{32}\neq 0,$ $  a_{11}= \frac{a_{31}a_{12}}{a_{32}} $
\\  \hline
$ \mathcal{A}_{10} $: $ e_{1}\cdot e_{3}= e_{1} $ & $ \left(
\begin{array}{ccc}
a_{11} & a_{12} & 0 \\
a_{21} & a_{22} & 0 \\
0 & 0 & 0 \\  
\end{array}
\right)
$; 
 $ \left(
\begin{array}{ccc}
0 & 0 & 0 \\
0 & 0 & 0 \\
a_{31} & 0 & 0 \\  
\end{array}
\right)a_{31}\neq 0
$; $ \left(
\begin{array}{ccc}
0 & 0 & 0 \\
0 & 0 & 0 \\
a_{31} & a_{32} & 0 \\  
\end{array}
\right) a_{32}\neq 0
$  \cr
 $ e_{2}\cdot e_{3}= e_{2} $ & 
  $ \left(
\begin{array}{ccc}
0 & 0 & a_{13} \\
a_{21} & 0 & 0 \\
0 & 0 & 0 \\  
\end{array}
\right)a_{13} \neq 0
 $; $ \left(
\begin{array}{ccc}
0 & 0 & a_{13} \\
0 & 0 & a_{23} \\
0 & 0 & 0 \\  
\end{array}
\right)
 a_{23}\neq 0$ \cr
$ e_{3}\cdot e_{1}= e_{1} $ & \cr
$ e_{3}\cdot e_{2}= e_{2} $ &\cr
$ e_{3}\cdot e_{3}= e_{3} $ &
\\ \hline
$ \mathcal{A}_{11} $: $ e_{1}\cdot e_{3}= e_{2} $ & $ \left(
\begin{array}{ccc}
a_{11} & a_{12} & 0 \\
a_{21} & a_{22} & 0 \\
0 & 0 & 0 \\  
\end{array}
\right) 
 $;
  $ \left(
\begin{array}{ccc}
0 & a_{12} & 0 \\
0 & 0 & 0 \\
0 & a_{32} & 0 \\  
\end{array}
\right) a_{32}\neq 0 
 $  \cr
$ e_{2}\cdot e_{3}= e_{2} $ & 
 $ \left(
\begin{array}{ccc}
0 & 0 & 0 \\
0 & 0 & a_{23} \\
0 & 0 & 0 \\  
\end{array}
\right) a_{23}\neq 0$;$ \left(
\begin{array}{ccc}
a_{11} & 0 & 0 \\
a_{21} & 0 & 0 \\
a_{31} & 0 & 0 \\  
\end{array}
\right)  a_{31}\neq 0 $  \cr
$ e_{3}\cdot e_{1}= e_{2} $ &$ \left(
\begin{array}{ccc}
a_{11} & 0 & 0 \\
a_{21} & 0 & 0 \\
a_{31} & a_{32} & 0 \\  
\end{array}
\right) a_{31} \neq 0, a_{32} \neq $ $ a_{11}= -a_{21} $ \cr
$ e_{3}\cdot e_{2}= e_{2} $ &$ \left(
\begin{array}{ccc}
0 & 0 & 0 \\
a_{21} & a_{22} & a_{23} \\
0 & 0 & 0 \\  
\end{array}
\right)a_{21} = -a_{22}\neq 0, a_{23}\neq 0 $  \cr
$ e_{3}\cdot e_{3}= e_{3} $ &
\\ \hline 
$ \mathcal{A}_{12} $: $ e_{1}\cdot e_{1}= e_{2} $ & $ \left(
\begin{array}{ccc}
0 & a_{12} & 0 \\
a_{21} & a_{22} & 0 \\
0 & 0 & 0 \\  
\end{array}
\right) a_{12}^{2} + a_{21}^{2}= a_{12}a_{21} $ \cr
$ e_{1}\cdot e_{3}= e_{1} $ &$ \left(
\begin{array}{ccc}
0& 0 & 0\\
0 & 0 &  \\
0& a_{32} & 0 \\  
\end{array}
\right) a_{32}\neq 0 $\cr
$ e_{2}\cdot e_{3}= e_{2} $ &\cr
$ e_{3}\cdot e_{1}= e_{1} $ &\cr
$ e_{3}\cdot e_{2}= e_{2} $ &\cr
 $ e_{3}\cdot e_{3}= e_{3} $ &
\\  \hline
\end{tabular}
\end{center}

\begin{center}
{\bf Table 8}: 
 Antisymmetric solutions and Frobenius algebra structures.
\begin{tabular}{c c c c}
\hline
Associative & Antisymmetric  & Frobenius algebras \cr
 algebra  & solutions     &  structures 
\\ \hline
  $ \mathcal{A}_{3} $ & $ a_{23}e_{2}\otimes e_{3} +  $ & $ e_{1}\ast e_{1}= e_{2}; e_{1}\ast e_{2}= e_{3}; e_{2}\ast e_{1}= e_{3}$; $ e_{1}\ast e^{\ast}_{2}= e^{\ast}_{1}$; $ e^{\ast}_{3}\ast e_{1}=  e^{\ast}_{2}$ \cr 
  &  $  a_{32}e_{3}\otimes e_{2} $  & $e^{\ast}_{3}\ast e^{\ast}_{3}= -2 a_{23}e^{\ast}_{1}$; $ e_{2}\ast e^{\ast}_{3}= e^{\ast}_{1}$; $ e^{\ast}_{3}\ast e_{2}=  e^{\ast}_{1}  $ \cr
 &  $  a_{32}= -a_{23} $  & $ e_{1}\ast e^{\ast}_{3}=  e^{\ast}_{2} -2a_{23}e_{3} $;
                        $  e^{\ast}_{2}\ast e_{1}= e^{\ast}_{1} -2a_{23}e_{3}  $
  \\ \hline
  $ \mathcal{A}_{4} $ &  $  a_{12}e_{1}\otimes e_{2} + $  &  $ e_{1}\ast e_{3}= e_{2}$; $ e_{2}\ast e_{3}= e_{2}$; $ e_{3}\ast e_{3}= e_{3}$; $ e^{\ast}_{3}\ast e_{3}= e^{\ast}_{3}$  \cr
&$ a_{21}e_{2}\otimes e_{1} +$&  $ e_{2}\ast e^{\ast}_{2}= a_{23}e_{2}; e_{1}\ast e^{\ast}_{2}= a_{23}e_{2} $;  $ e^{\ast}_{1}\ast e_{3}= a_{12}e_{2} $ \cr
                   & $ a_{23}e_{2}\otimes e_{3} +$ & $ e^{\ast}_{2}\ast e^{\ast}_{1}= a_{12}e^{\ast}_{3};$ $  e^{\ast}_{3}\ast e^{\ast}_{2}= a_{23} e^{\ast}_{3} $ \cr
                   &$  a_{32}e_{3}\otimes e_{2}  $ &$ e^{\ast}_{2}\ast e_{3}= -a_{12}e^{\ast}_{2} + a_{23}e_{3}$; $ e_{3}\ast e^{\ast}_{3}= e^{\ast}_{3} + a_{23}e_{2} $ \cr
                        & $ a_{21}= -a_{12} $ & $ e^{\ast}_{2}\ast e_{1}= e^{\ast}_{3} + a_{23}e_{2};$  $ e_{3}\ast e^{\ast}_{2}= e^{\ast}_{1} + e^{\ast}_{2} -a_{12}e_{2} + a_{12}e_{1}$  \cr
                        & $ a_{32}= -a_{23} $ &$ e^{\ast}_{2}\ast e_{2}= e^{\ast}_{3} + a_{23}e_{2} $; $ e^{\ast}_{2}\ast e^{\ast}_{2}= a_{23}e^{\ast}_{1} + a_{23}e^{\ast}_{2} -a_{12}e^{\ast}_{3} $
\\ \hline
$ \mathcal{A}_{5} $ &  $  a_{12}e_{1}\otimes e_{2} + $  &  $ e_{2}\ast e_{3}= e_{2}$; $ e_{3}\ast e_{1}= e_{1}$; $ e_{3}\ast e^{\ast}_{3}= e^{\ast}_{3}$; $ e_{3}\ast e^{\ast}_{2}= e^{\ast}_{2}$;  \cr
                         &  $ a_{21}e_{2}\otimes e_{1}+$  & $ e^{\ast}_{1}\ast e_{3}= e^{\ast}_{1}; e_{3}\ast e^{\ast}_{3}= e^{\ast}_{3}$; $ e^{\ast}_{1}\ast e_{1}= a_{13}e_{1} $;  \cr
                   & $ a_{13}e_{1}\otimes e_{3} + $ & $e_{2}\ast e^{\ast}_{1}= a_{13}e_{2}; e_{3}\ast e^{\ast}_{1}= a_{13}e_{3}; e^{\ast}_{1}\ast e^{\ast}_{1}= a_{13} e^{\ast}_{1};$ \cr
                   &$ a_{31}e_{3}\otimes e_{1}  $ & $e^{\ast}_{1}\ast e^{\ast}_{2}= a_{13} e^{\ast}_{2};$ $ e^{\ast}_{2}\ast e_{2}=  e^{\ast}_{3} + a_{13}e_{1};$ $e^{\ast}_{1}\ast e^{\ast}_{3}= a_{13} e^{\ast}_{3};$  \cr
                        & $ a_{21}= -a_{12} $ & $ e_{1}\ast e^{\ast}_{1}=  e^{\ast}_{3} + a_{13}e_{1}$; $ e^{\ast}_{3}\ast e_{3}=  e^{\ast}_{3} + a_{13}e_{1}.$   \cr
 &  $ a_{13}= -a_{31} $ & \cr
\cline{2-3}  &  $a_{12}e_{1}\otimes e_{2}+$  &  $e_{2}\ast e_{3}= e_{2}$; $ e_{3}\ast e_{1}= e_{1}$; $e_{3}\ast e^{\ast}_{3}= e^{\ast}_{3}$  \cr
                         &  $ a_{21}e_{2}\otimes e_{1} + $  &  $e_{2}\ast e^{\ast}_{2}= a_{23}e^{\ast}_{2}$; $ e^{\ast}_{2}\ast e_{1}= a_{23}e_{1}; e^{\ast}_{2}\ast e_{3}= a_{23}e_{3};$\cr
                   & $ a_{23}e_{2}\otimes e_{3} + $ & $ e^{\ast}_{1}\ast e^{\ast}_{2}= a_{23}e^{\ast}_{1}; e^{\ast}_{2}\ast e^{\ast}_{2}= a_{23}e^{\ast}_{2}; e^{\ast}_{2}\ast e^{\ast}_{3}= a_{23}e^{\ast}_{3};$ \cr
                   &$ a_{32}e_{3}\otimes e_{2}$ &  $ e^{\ast}_{3}\ast e^{\ast}_{2}= a_{23} e^{\ast}_{3};$ $ e^{\ast}_{2}\ast e_{2}=  e^{\ast}_{3} + a_{23}e_{2}$;  \cr
                        & $ a_{21}= -a_{12} $ &$ e^{\ast}_{1}\ast e_{3}=  e^{\ast}_{1} - a_{12}e_{2}   $; $ e_{1}\ast e^{\ast}_{1}=  e^{\ast}_{3} + a_{23}e_{2};$ \cr
 &  $ a_{23}= -a_{32} $ &  $ e^{\ast}_{3}\ast e_{3}=  e^{\ast}_{3} + a_{23}e_{2};$ $ e_{3}\ast e^{\ast}_{3}=  e^{\ast}_{3} + a_{23}e_{2};$\cr
 &&$e^{\ast}_{2}\ast e^{\ast}_{1}= -a_{12} e^{\ast}_{3}$; $ e_{3}\ast e^{\ast}_{2}=  e^{\ast}_{2} + a_{23}e_{3} -a_{12}e_{1}$ .
\\ \hline
$ \mathcal{A}_{6} $ & $ a_{23}e_{2}\otimes e_{3} +  $ & $ e_{3}\ast e_{1}= e_{2}; e_{3}\ast e_{2}= e_{2};$ $ e_{3}\ast e_{3}= e_{3}; e^{\ast}_{3}\ast e_{3}= e^{\ast}_{3}$; \cr 
  &  $  a_{32}e_{2}\otimes e_{1} $  & $ e_{3}\ast e^{\ast}_{3}= e^{\ast}_{3};  e^{\ast}_{2}\ast e_{1}= a_{23}e_{2}; e^{\ast}_{2}\ast e_{2}= a_{23}e_{2};$ \cr
  & $ a_{32} = - a_{23} $   &$   e^{\ast}_{2}\ast e^{\ast}_{3}= a_{23}e^{\ast}_{3}$; $ e^{\ast}_{2}\ast e_{3}=  e^{\ast}_{1} + e^{\ast}_{2} $  \cr
 && $ e_{1}\ast e^{\ast}_{2}=  e^{\ast}_{3} + a_{23}e_{2};$  $ e_{2}\ast e^{\ast}_{2}=  e^{\ast}_{3} + a_{23}e_{2};$ \cr
 && $ e^{\ast}_{2}\ast e^{\ast}_{2}=  a_{23}e^{\ast}_{1} + a_{23}e^{\ast}_{2}.$
    \\ \hline
$ \mathcal{A}_{7} $ & $ a_{12}e_{1}\otimes e_{2} +  $ & $ e_{1}\ast e_{2}= e_{1}; e_{2}\ast e_{2}= e_{2};$ $ e_{3}\ast e_{1}= e_{1}; e_{3}\ast e_{3}= e_{3}$ \cr 
  &  $  a_{21}e_{2}\otimes e_{1}$  &$ e^{\ast}_{1}\ast e^{\ast}_{1}= a_{12}e^{\ast}_{1};  e_{3}\ast e^{\ast}_{2}= -a_{12}e_{1}; e^{\ast}_{1}\ast e_{2}= a_{12}e_{2}; e^{\ast}_{2}\ast e_{2}=e^{\ast}_{2}$  \cr
  & $ a_{12} = - a_{21} $ & $ e^{\ast}_{1}\ast e_{1}=  e^{\ast}_{2} + a_{12}e_{1};  e_{3}\ast e^{\ast}_{3}= e^{\ast}_{3}; e^{\ast}_{3}\ast e_{3}=e^{\ast}_{3}$  \cr
                   &  & $ e^{\ast}_{1}\ast e_{3}=  e^{\ast}_{1} - a_{12}e_{2};$ $ e_{1}\ast e^{\ast}_{1}=  e^{\ast}_{3} + a_{12}e_{1}   $ \cr
                        &  &   $ e^{\ast}_{2}\ast e^{\ast}_{1}= a_{12}e^{\ast}_{2} - a_{12}e^{\ast}_{3} + e^{\ast}_{1};$ $ e_{2}\ast e^{\ast}_{2}=  e^{\ast}_{2} + a_{12}e_{1}$ \cr
       \cline{2-3}
 &  $ a_{13}e_{1}\otimes e_{3} + $ & $ e_{1}\ast e_{2}= e_{1}; e_{2}\ast e_{2}=  e_{2};$ $ e_{3}\ast e_{1}=  e_{1}; e_{3}\ast e_{3}=  e_{3}  $\cr
 & $  a_{31}e_{3}\otimes e_{1}$ & $ e^{\ast}_{1}\ast e^{\ast}_{1}= a_{13}e^{\ast}_{1}; e^{\ast}_{3}\ast e_{2}=- a_{13}e_{1}; e_{3}\ast e^{\ast}_{1}= a_{13}e_{3} $ \cr
& $ a_{31}= -a_{13} $ & $e_{3}\ast e^{\ast}_{3}= e^{\ast}_{3};$ $ e^{\ast}_{1}\ast e_{3}= e^{\ast}_{1}; e^{\ast}_{2}\ast e_{2}= e^{\ast}_{2};$ $ e_{2}\ast e^{\ast}_{2}= e^{\ast}_{2}$\cr
                   &  & $ e^{\ast}_{1}\ast e_{1}= e^{\ast}_{1}+ a_{13}e_{1};$ $ e_{2}\ast e^{\ast}_{1}= e^{\ast}_{1}- a_{13}e_{3}  $ \cr
                        &  & $  e_{1}\ast e^{\ast}_{1}= e^{\ast}_{3}+ a_{13} e^{\ast}_{1};$ $ e^{\ast}_{1}\ast e^{\ast}_{3}= -a_{13} e^{\ast}_{2} + a_{13} e^{\ast}_{3} $
 \\ \hline
$ \mathcal{A}_{8} $ &  $ a_{12}e_{1}\otimes e_{3} + $  & $ e_{1}\ast e_{3}=  e_{1}; e_{2}\ast e_{3}=  e_{2}$; $ e_{3}\ast e_{1}=  e_{1}; e_{3}\ast e_{3}=  e_{3}   $ \cr
& $ a_{21}e_{2}\otimes e_{1}  $ & $ e^{\ast}_{1}\ast e_{1}=  e^{\ast}_{3}; e^{\ast}_{2}\ast e_{2}=  e^{\ast}_{3} $; $ e_{3}\ast e^{\ast}_{2}=  e^{\ast}_{2}; e^{\ast}_{1}\ast e_{3}=  e^{\ast}_{1}   $ \cr
& $ a_{21}= -a_{12} $ & $ e_{1}\ast e^{\ast}_{1}=  e^{\ast}_{3}; e^{\ast}_{3}\ast e_{3}=  e^{\ast}_{3}; e_{3}\ast e^{\ast}_{3}= e^{\ast}_{3};$ \cr
&& $ e^{\ast}_{1}\ast e^{\ast}_{2}= -a_{12}e^{\ast}_{3};$ $ e^{\ast}_{2}\ast e_{3}= -a_{12}e_{1};$ \cr
&& $e_{3}\ast e^{\ast}_{1}=  e^{\ast}_{1} -a_{12}e_{2}$ \cr
\cline{2-3}
 &  $ a_{23}e_{2}\otimes e_{3} + $  & $ e_{1}\ast e_{3}=  e_{1}; e_{2}\ast e_{3}=  e_{2}$; $ e_{3}\ast e_{1}=  e_{1}; e_{3}\ast e_{3}=  e_{3}   $ \cr
& $ a_{32}e_{3}\otimes e_{2}$ & $ e_{3}\ast e^{\ast}_{1}=  e^{\ast}_{1}; e_{3}\ast e^{\ast}_{2}=  e^{\ast}_{2}$; $ e^{\ast}_{1}\ast e_{3}=  e^{\ast}_{1}; e^{\ast}_{3}\ast e_{3}=  e^{\ast}_{3}   $ \cr
& $ a_{32}= -a_{23} $ & $ e^{\ast}_{3}\ast e^{\ast}_{2}= a_{23}e^{\ast}_{3}; e^{\ast}_{2}\ast e^{\ast}_{2}= a_{23}e^{\ast}_{2}  $ \cr
& & $ e^{\ast}_{1}\ast e^{\ast}_{2}= a_{23}e^{\ast}_{1}; e^{\ast}_{2}\ast e^{\ast}_{1}= a_{23}e^{\ast}_{1}; e^{\ast}_{2}\ast e_{3}= a_{23}e_{3}$\cr
&& $ e_{2}\ast e^{\ast}_{2}= a_{23}e_{2};$ $e^{\ast}_{2}\ast e_{1}= a_{23}e_{1}; e_{1}\ast e^{\ast}_{2}= a_{23}e_{1}$\cr
&& $ e_{3}\ast e^{\ast}_{3}= e^{\ast}_{3} + a_{23}e^{\ast}_{2};$ $ e_{1}\ast e^{\ast}_{1}= e^{\ast}_{3} + a_{23}e_{1} $ \cr
&& $ e^{\ast}_{2}\ast e_{2}= e^{\ast}_{3} + a_{23}e_{2};$ $ e^{\ast}_{1}\ast e_{1}= e^{\ast}_{3} + a_{23}e_{2} $ 
  \\ \hline
\end{tabular}
\end{center}

 \begin{center}
\begin{tabular}{c c c c}
  \\ \hline
  $ \mathcal{A}_{9} $ & $ a_{12}e_{1}\otimes e_{2} +  $ & $ e_{2}\ast e_{3}= e_{2}; e_{3}\ast e_{1}= e_{1};$ $ e_{3}\ast e_{2}= e_{2}; e_{3}\ast e_{3}= e_{3}$ \cr 
  &  $  a_{21}e_{2}\otimes e_{1} $  & $ e_{3}\ast e^{\ast}_{2}= e^{\ast}_{2};  e^{\ast}_{2}\ast e_{2}= e^{\ast}_{3};$ $ e_{1}\ast e^{\ast}_{1}= e^{\ast}_{3};  e^{\ast}_{1}\ast e_{3}= e^{\ast}_{1}  $  \cr
  & $ a_{12} = - a_{21} $   &  $ e_{2}\ast e^{\ast}_{2}=e^{\ast}_{3}; e^{\ast}_{3}\ast e_{3}=e^{\ast}_{3};$ $ e_{3}\ast e^{\ast}_{3}=  e^{\ast}_{3}; e_{3}\ast e^{\ast}_{1}= a_{12}e_{2}  $  \cr
                         &    & $ e^{\ast}_{1}\ast e^{\ast}_{2}=  a_{12}e^{\ast}_{3};$ $ e^{\ast}_{2}\ast e_{3}=  e^{\ast}_{2} + a_{12}e_{1}   $ \cr
 \cline{2-3}
  &  $ a_{13}e_{1}\otimes e_{3} + $ & $e_{2}\ast e_{3}= e_{2}; e_{3}\ast e_{1}=  e_{1}; e_{3}\ast e_{2}=  e_{2}; e_{3}\ast e_{3}=  e_{3}$\cr
 & $  a_{31}e_{3}\otimes e_{1}   $ & $ e_{3}\ast e^{\ast}_{2}= e^{\ast}_{2}; e^{\ast}_{2}\ast e_{3}= e^{\ast}_{2};$ $ e_{3}\ast e^{\ast}_{3}= e^{\ast}_{3}; e^{\ast}_{1}\ast e_{1}= a_{13}e_{1}$ \cr
& $ a_{31}= -a_{13} $ & $ e_{3}\ast e^{\ast}_{1}= a_{11}e_{1}; e_{2}\ast e^{\ast}_{1}= a_{13}e_{2};$ $ e^{\ast}_{1}\ast e_{2}= a_{13}e_{2}; e_{3}\ast e^{\ast}_{1}= a_{13}e_{3}$ \cr
                        &    &  $ e^{\ast}_{1}\ast e^{\ast}_{2}= a_{13}e^{\ast}_{2}; e^{\ast}_{2}\ast e^{\ast}_{1}= a_{13}e^{\ast}_{2};$ $ e^{\ast}_{1}\ast e^{\ast}_{3}= a_{13}e^{\ast}_{3}$  \cr
                   &  & $ e^{\ast}_{1}\ast e^{\ast}_{1}= a_{13}e^{\ast}_{1}+ a_{11}e^{\ast}_{3};$  $ e^{\ast}_{2}\ast e_{2}= e^{\ast}_{3} + a_{13}e_{1}  $ \cr
                        &  & $  e_{1}\ast e^{\ast}_{1}= e^{\ast}_{3}+ a_{13} e_{1};$  $  e^{\ast}_{1}\ast e_{3}= e^{\ast}_{1}+ a_{11} e_{1}    $ \cr
 &   & $ e_{2}\ast e^{\ast}_{2}= e^{\ast}_{3} + a_{13} e_{1};$ 
 $ e^{\ast}_{3}\ast e_{3}= e^{\ast}_{3} + a_{13} e_{1} $ 
\\ \hline 
$ \mathcal{A}_{10} $ &  $ a_{12}e_{1}\otimes e_{2} + $  & $ e_{1}\ast e_{3}=  e_{1}; e_{2}\ast e_{3}=  e_{2};$ $ e_{3}\ast e_{1}=  e_{1}; e_{3}\ast e_{2}=  e_{2}$  \cr
&  $ a_{21}e_{2}\otimes e_{1}  $  & $ e_{3}\ast e_{3}=  e^{\ast}_{3}; e_{3}\ast e^{\ast}_{1}=  e^{\ast}_{1};$ $ e^{\ast}_{1}\ast e_{1}=  e^{\ast}_{3}; e^{\ast}_{2}\ast e_{2}=  e^{\ast}_{3}   $ \cr
& $ a_{21}= -a_{12} $ & $ e_{3}\ast e^{\ast}_{2}=  e^{\ast}_{2}; e^{\ast}_{1}\ast e_{3}=  e^{\ast}_{1};$ $ e_{1}\ast e^{\ast}_{1}= e^{\ast}_{3}; e^{\ast}_{2}\ast e_{3}= e^{\ast}_{2}   $\cr
&& $ e_{2}\ast e^{\ast}_{2}= e^{\ast}_{3}; e_{3}\ast e^{\ast}_{3}= e^{\ast}_{3};$ $e^{\ast}_{3}\ast e_{3}=  e^{\ast}_{3} $ 
\\ \hline

$ \mathcal{A}_{11} $ &  $ a_{12}e_{1}\otimes e_{2} + $  & $ e_{1}\ast e_{3}=  e_{2}; e_{2}\ast e_{3}=  e_{2};$ $ e_{3}\ast e_{1}=  e_{2}; e_{3}\ast e_{2}=  e_{2}   $ \cr
&  $ a_{21}e_{2}\otimes e_{1}$  & $ e_{3}\ast e_{3}=  e_{3}; e^{\ast}_{2}\ast e_{1}=  e^{\ast}_{3};$ $ e^{\ast}_{2}\ast e_{2}=  e^{\ast}_{3}; e_{1}\ast e^{\ast}_{2}=  e^{\ast}_{3}   $ \cr
& $ a_{21}= -a_{12} $ &  $ e_{2}\ast e^{\ast}_{2}= e^{\ast}_{3}; e^{\ast}_{3}\ast e_{3}= e^{\ast}_{3};$  $ e_{3}\ast e^{\ast}_{3}= e^{\ast}_{3}; e_{3}\ast e^{\ast}_{1}= a_{12}e_{2} $\cr
&& $e^{\ast}_{1}\ast e_{3}= a_{12}e_{2}; e^{\ast}_{1}\ast e^{\ast}_{2}= a_{12}e^{\ast}_{3}; e^{\ast}_{2}\ast e^{\ast}_{2}= -2a_{12}e^{\ast}_{3} $ \cr
&& $ e^{\ast}_{2}\ast e^{\ast}_{1}= a_{12}e^{\ast}_{3};$ $ e_{3}\ast e^{\ast}_{2}= e^{\ast}_{1} + a_{12}e_{1} + e^{\ast}_{2}-2a_{12}e_{2} $ \cr
&& $ e^{\ast}_{2}\ast e_{3}= e^{\ast}_{2} -2a_{12}e_{2} + e^{\ast}_{1} +a_{12}e_{1} $ 
  \\ \hline 
$ \mathcal{A}_{12} $ & $ r=0 $ & $ e_{1}\ast e_{1}= e_{2}; e_{1}\ast e_{3}= e_{1};$  $ e_{2}\ast e_{3}= e_{2}; e_{3}\ast e_{1}= e_{1}$ \cr
  &    & $ e_{3}\ast e_{2}= e_{2};  e_{3}\ast e_{3}= e_{3};$  $ e_{1}\ast e^{\ast}_{2}= e^{\ast}_{1};  e_{3}\ast e^{\ast}_{1}= e^{\ast}_{1}  $ \cr
  &  & $ e_{3}\ast e^{\ast}_{2}= e^{\ast}_{2}; e_{1}\ast e^{\ast}_{1}=e^{\ast}_{3};$      $ e_{2}\ast e^{\ast}_{2}= e^{\ast}_{3}; e_{3}\ast e^{\ast}_{3}= e^{\ast}_{3}  $  \cr
                         &    & $ e^{\ast}_{1}\ast e_{1}=  e^{\ast}_{3}; e^{\ast}_{2}\ast e_{1}=  e^{\ast}_{1};$ $ e^{\ast}_{2}\ast e_{1}=  e^{\ast}_{3}; e^{\ast}_{1}\ast e_{3}=  e^{\ast}_{1}  $\cr
                     &    & $ e^{\ast}_{2}\ast e_{3}=  e^{\ast}_{2}; e^{\ast}_{3}\ast e_{3}=  e^{\ast}_{3}  $

  \\ \hline 

\end{tabular}
\end{center}

 \section{Concluding remarks}
 In this work, we have considered a three-dimensional associative algebra $\mathcal{H}$ consisting of the $3\times 3$ strictly upper triangular matrices whose the commutator is the Heisenberg Lie algebra. 
Then, we have investigated the concrete case  of extended associative Heisenberg algebra ${\mathcal H}_{\alpha},$ realized double constructions of its  Frobenius like algebras, Connes cocycles and their antisymmetric infinitesimal bialgebra stuctures. We have 
explicitly determined  an associative algebra structure on the direct sum 
$ \mathcal{H}_{\alpha} \oplus \mathcal{H}^{\ast}_{\alpha}  $ of the underlying vector spaces of $ \mathcal{H}_{\alpha} $ and $ \mathcal{H}^{\ast}_{\alpha} $ such 
that $ (\mathcal{H}_{\alpha}, \cdot) $ and ($ \mathcal{H}^{\ast}_{\alpha}, \circ $) are subalgebras and the symmetric bilinear form  $ \mathcal{B}(\cdot, \cdot ) $ defined on 
$ \mathcal{H}_{\alpha}\oplus \mathcal{H}^{\ast}_{\alpha}$  
is invariant. Hence,  $ (\mathcal{H}_{\alpha}\oplus \mathcal{H}^{\ast}_{\alpha}, \mathcal{B}) $ 
has been proved to be  a symmetric Frobenius algebra.
Besides, we have explicitly constructed a matched pair of associative algebras 
$ \mathcal{H}_{\alpha} $ and $ \mathcal{H}^{\ast}_{\alpha} $ and the associated double structure 
on the direct sum of the associative algebras $ \mathcal{H}_{0} $ and $ \mathcal{H}^{\ast}_{0} $
 in the matched pair. Moreover, we have established a dendriform algebra structure on $ \mathcal{H}^{\ast}_{0} $.
Finally we have built a double construction of the Connes cocycle associated to 
associative algebras  $(\mathcal{H}_{0}, \ast_{\mathcal{H}_{0}})$ and $(\mathcal{H}^{\ast}_{0}, \ast_{\mathcal{H}^{\ast}_{0}})$.
The realization of the  Heisenberg associative double associated to the constructed 
Frobenius algebras, as well as the study of its relations with the  antisymmetric 
solutions  of the associative  Yang-Baxter equation have been performed. Besides, has 
followed  a classification of the solutions of $3-$dimensional non decomposable associative  Yang-Baxter 
equations. We have also provided   important examples of compatible dendriform 
algebras on the extended associative Heisenberg algebra $\mathcal{H}_{\alpha}$. Besides, 
we have proceeded to do the classification of solutions of  $D-$equations in the considered  
dendriform algebras, and to give the corresponding double constructions of Connes 
cocycles.

 \section*{Aknowledgement}
 This work is partially supported by the Abdus Salam International Centre for Theoretical
 Physics (ICTP, Trieste, Italy) through the Office of External Activities (OEA) - Prj-15. The
 ICMPA is also in partnership with the Daniel Iagolnitzer Foundation (DIF), France.

\end{document}